\documentclass[a4paper]{article}
\usepackage{RR}
\usepackage{hyperref}
\RRdate{Mars 2009}
\usepackage{amsmath,amssymb,epsf,graphicx,a4wide}
\RRetitle{Estimating nonlinearities in twophase flow in porous media}
\RRtitle{Estimation des non-lin\'earit\'es pour des \'ecoulements diphasiques en milieu poreux}
\RRauthor{Jianfeng Zhang\thanks{INRIA, Domaine de Voluceau, 78153 Le Chesnay Cedex, France}
\and Guy Chavent\thanks{CEREMADE, Universit\'e Paris-Dauphine, France}\footnotemark[1] \and J\'er\^ome Jaffr\'e\footnotemark[1]}
\authorhead{J. Zhang, G. Chavent \& J. Jaffr\'e}

\newcommand{\R}{\mathbb R}

\newcommand{\p}{\partial}

\newtheorem{theo}{Theorem}
\newtheorem{definition}{Definition}
\RRNo{6892}
\RRdate{Mars 2009}
\RRabstract{
In order to analyze numerically inverse problems several techniques
 based on linear and nonlinear stability analysis are presented. These
techniques are illustrated on the problem of estimating mobilities and
capillary pressure in one-dimensional two-phase displacements in
porous media that are performed in laboratories. This is an example of
the problem of estimating nonlinear coefficients in a system of
nonlinear partial differential equations.}
\RRkeyword{ flow in porous media, inverse problem, estimation of
nonlinear coefficients.}
\RRresume{Afin d'analyser num\'eriquement des probl\`emes inverses on
pr\'esente plusieurs techniques bas\'ees sur l'analyse de stabilit\'e lin\'eaire et 
non-lin\'eaire. Ces techniques sont pr\'esent\'ees pour le probl\`eme d'estimation des mobilit\'es et de la pression capillaire dans des d\'eplacements diphasiques unidimensionnels en milieu poreux
r\'ealis\'es en laboratoire. C'est un exemple de probl\`eme d'estimation des coefficients non-lin\'eaires dans un
syst\`eme d'\'equations aux d\'eriv\'ees partielles non-lin\'eaires.}
\RRmotcle{\'ecoulement en milieu poreux, probl\`eme inverse, estimation des coefficients non-lin\'eaires}
\RRprojets{Estime}
\RRtheme{\THNum}
\RCParis

\begin{document}
\RRNo{6892}
\makeRR

\section{Introduction}
Multiphase flow in porous media is modelled by a set of nonlinear
partial differential equations
equations and it provides a very good practical example for the
inverse problem of estimating
nonlinear coefficients in nonlinear partial differential
equations. The standard problem in petroleum engineering is to
estimate the relative permeabilities and capillary pressure curves from
laboratory experiments which consists of displacing a resident phase
by injecting the other
\cite{chalem74,chacoh78,chacohespy80,watgavsei84,kerwat86,chachajafliu90,chachajafliubou92,naemanbrusnor00}.
The relative permeabilities and the capillary pressure are functions
of the saturation of one of the phases.
More recently experiments where the displacement is due to
centrifugation were designed in order to improve the estimation of the
capillary pressure function
\cite{chaforzhachalen92,vignes93,normejyanwat93,chachazha94}. 
Three-phase flow were also considered in
\cite{chajafjeg99,mejwatnor96}. In this case the relative
permeabilities and the capillary pressure are functions of two variables.
In hydrogeology the Richards equation is often used and the
problem of estimating its coefficients is considered in
\cite{iglkna97,bitkna}.

Without trying to give a complete review we can add to this
bibliography several interesting contributions
\cite{duc93,tai95,griman00,grikolmannor01} and two reviews for
parameter estimation in multiphase flow \cite{ewipilwadwat95,sun}.

In this paper we present several ingredients for a successful
numerical estimation of the relative permeabilities and the capillary
pressure. In Section 2 we introduce the mathematical model for
two-phase flow and in Section 3 we set the parameter estimation
problem as a minimization problem. Multiscale parameterization is
adressed in Section 4. Some linear analysis of the problem is
presented in Section 5 and confidence intervals are calculated in
Section 6 using edgehog extremal solutions. Techniques for nonlinear
analysis are presented in Section 7 and implemented numerically in
Section 8.

\section{A model for a two-phase displacement in porous media}
\begin{figure}[htbp]
 \begin{center}
\begin{picture}(400,80)(0,0)
    \includegraphics[height=3cm]{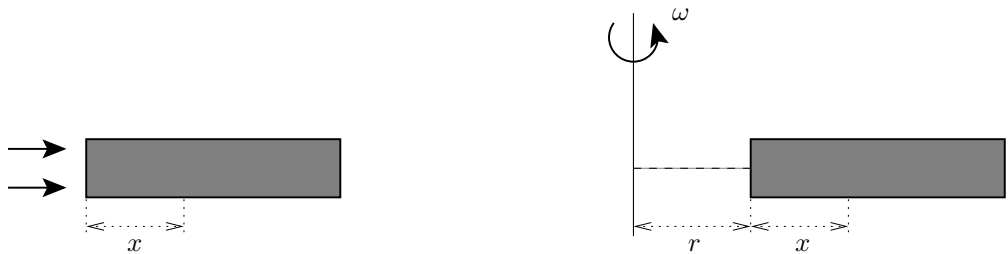}
\put(-120,-5){$r$}
\put(-80,-5){$x$}
\put(-330,-5){$x$}
\put(-126,82){$\omega$}
\end{picture}
 \end{center}
  \caption{Two-phase displacement by injection (right) or by
    centrifugation (left)} 
  \label{fig1}
\end{figure}

In several laboratories core samples collected from oil fields are
analyzed to determine their flow properties. A typical experiment
consist in displacing a resident wetting fluid (subscript $w$), say
water, by a nonwetting fluid (subscript $nw$), say oil. The displacement may be
driven by injecting the nonwetting fluid through one extremity of
the core  or by centrifugal forces. These experiments are sketched in
Fig. \ref{fig1}.

Two-phase displacement is governed by a generalized Darcy's law and,
in laboratory experiments, it is usually
assumed to be incompressible and one-dimensional. Using the global
pressure formulation~\cite{jaf1} the displacement is modelled by the
following nonlinear equation :
\begin{equation} \begin{array}{l}
\displaystyle{\phi\frac{\p S}{\p t} + \frac{\p q_w}{\p x}} = 0,\\
q_w = - Ka(S)\displaystyle{\frac{\p S}{\p x}}+q_T(t) b_T(S) + q_G b_G(S),
\end{array}
\label{eqsaturation}
\end{equation}
where $S=S_w ((0 \leq S \leq 1)$ is the saturation of the wetting
fluid, $q_w$ its Darcy velocity, $\phi$ is the
porosity of the rock and $K$ is its absolute permeability. The total
flow rate $q_T = q_w + q_{nw}$ is the sum of the flow rate of
the two phases and is given by
\begin{equation}
q_T = -Kd(S)[\frac{\p P}{\p x}-\rho(S)H].
\label{eqpression}
\end{equation}
The global pressure $P$ \cite{jaf1} is given by
\begin{equation}
\label{globalpressureeq}
  P = \frac{1}{2}(p_w + p_{nw}) + \gamma(S), 
\end{equation}
with $p_w, p_nw$ the phase pressures and $\gamma$ defined below.

$H$ is a gravity or centrifugation function. In case of a standard displacement
$H(x) = g \nabla Z(x)$ where g is the Newton constant and $Z$ is the depth at
the location $x$. In case of a displacement by centrifugation $H(x) = \omega^2(r +x)$
where $\omega$ is the angular speed, $r$ is 
the distance from the rotation axis to the
closest extremity of the core (see Fig. \ref{fig1}).
 The total flow rate $q_T$
is independent of $x$ because of the incompressibility assumption.

The gravity or centrifugation field $q_G$ is given by
\begin{equation}
q_G = K \displaystyle{\frac{\rho_w+\rho_{nw}}{2}}\, H
\end{equation}
where $\rho_w$ and $\rho_{nw}$ are the densities of the wetting fluid and
the nonwetting fluid respectively.

We have introduced above the
coefficients $a$, $b_T$, $b_G$, $d$, $\rho, \gamma$ which are functions of the
saturation $S$. They relate to the relative permeability functions $kr_w$ and $kr_{nw}$
and to the capillary pressure function $p_c = p_w - p_{nw}$ through
the following relations :
\begin{displaymath}
\begin{array}{lll}
a={\displaystyle \frac{k_wk_{nw}}{k_w+k_{nw}}}p'_c, 
&b_T={\displaystyle \frac{k_w}{k_w+k_{nw}}}, &
\rule{0cm}{1cm} b_G={\displaystyle \frac{k_wk_{nw}}{k_w+k_{nw}}\;
\frac{\rho_w-\rho_{nw}}{{\frac{1}{2}}(\rho_w+\rho_{nw})}},\\
 d=k_w+k_{nw}, &
\rule{0cm}{1cm}\rho={\displaystyle
  \frac{k_w\rho_w+k_{nw}\rho_{nw}}{k_w+k_{nw}}}, & 
\gamma=\displaystyle{\int_0^S} (b_T(s)-\frac{1}{2}) \displaystyle{\frac{d 
p_c}{d S}}, \\
k_i = {\displaystyle \frac{kr_i}{\mu_i}}, \; i=w,nw.&
\end{array}
\end{displaymath}
where $\mu_i, i=w,nw$ are the viscosities of the two phases.

The relative permeabilities $kr_w$ and $kr_{nw}$ and the capillary
pressure $p_c$ are functions of the
saturation which satisfy the following physical
properties :
\begin{equation}
  \begin{array}{l}
kr_w \geq 0,\; kr_w \, \mbox{increasing},\;\;
kr_{nw} \geq 0,\; kr_{nw} \,\mbox{decreasing},\;\;
kr_w(0)=kr_{nw}(1)=0,\\ 
p_c \geq 0, \; p_c \,\mbox{decreasing}.
  \end{array}
  \label{properties}
\end{equation}

To equations (\ref{eqsaturation}), (\ref{eqpression}) we add various boundary
conditions depending on the experiments \cite{thesezhang}.

\section{The parameter estimation problem}
The problem is to estimate relative permeabilities and pressure capillary
functions. For this purpose experiments are set up so the following
measurements are available :
\begin{enumerate}
\item Local measurements : saturations $S^m_{k,i}$ are measured at
      different times $t_k$ and different locations $x_i$;
\item Global measurements : cumulated productions 
      $Q^m_k =   \int_0^{t_k} \phi_w(t) dt$  and pressure drops
      $\Delta P_k$ are measured at different times $t_k$.
\end{enumerate}
The problem of estimating the mobility and pressure capillary
functions is set as the problem of minimizing the error function $J$ :
\begin{eqnarray}
\lefteqn{J(kr_w,kr_{nw},p_c) = \sum_k\sum_i w_s^{k,i} (S^c_{k,i}-S^m_{k,i})^2 +
  \sum_k w_q^k (Q^c_k-Q^m_k)^2 } \hspace*{7.5cm} \nonumber \\ 
& & + \sum_k w_p^k (\Delta P^c_k-\Delta P^m_k)^2 .
  \label{eq3.1}
\end{eqnarray}
Here the superscripts $m$ and $c$ refer respectively to ``measured''
and  ``calculated'', and $w_s^{k,i}$, $w_q^k$ and $w_p^k$ are weights
given to the 
measurements. The function $J$ measures the difference between the
measured quantities and that calculated with the model using the parameters
$kr_w,kr_{nw},p_c$. 

The choice of parameterization of $kr_w,kr_{nw},p_c$ will be discussed
in the next section and is crucial for a successful estimation.
If computational costs are taken into account, the choice of the
parameterization determines also the choice of the minimization method.
If the choice of the parameterization gives a small number of
parameters, say smaller than 15, then Levenberg-Marquart or trust
region methods \cite{more78}
can be used and be cost efficient.
However, when the number of parameters becomes large, then these
optimization methods become too expensive and quasi-Newton methods \cite{bongillemsag03}
using the gradient of $J$ calculated with the adjoint method \cite{chavent74,sun} becomes
the right choice. The drawback of this method is the difficulty to
calculate the gradient and to implement this calculation. Therefore
techniques of automatic differentiation were developped
\cite{jegou97,gri00} and sophisticated software like Tapenade and
Adifor are now available.

Now we introduce some notations. We denote by ${\cal A}$ the set of
admissible parameters $p=(kr_w,kr_{nw},p_c)$, that is parameters that
satisfy properties (\ref{properties}). ${\cal A}$ is a subset
of a set ${\cal U}$. The direct mapping is the mapping $\varphi$ 
\begin{equation}
 \begin{array}[t]{cccc}
       \varphi:   &{\cal A} \subset \cal U          & \longmapsto & \cal O  \\ 
             & p =(kr_w,kr_{nw},p_c)  & \longmapsto & \varphi({p})= \varphi(u(p))=(S^c,Q^c,\Delta P^c)
 \end{array}
  \label{defphi}
\end{equation}
which maps ${\cal A}$ into the Hilbert space ${\cal O}$ of
observations. To solve the direct problem is to calculate
  $\varphi(p)$ for a given parameter $p$. This includes solving
  equations (\ref{eqsaturation}), (\ref{eqpression}) with appropriate
  initial and boundary conditions.

Let measurements $z = (Q^m,S^m,{\Delta P}^m)$ be given and let us write the
error function $J$ defined in (\ref{eq3.1}) in compact form as
\begin{equation}
      J({p})=\parallel \varphi ({p})- z {\parallel}^2_W ,
\label{cout}       
\end{equation} 
where $\parallel .{\parallel}_W $ is the weighted norm for $\cal O$
with $W$ the diagonal matrix of the weights given to the various measurements.
The inverse problem is to minimize $J$ over the set of admissible
parameters:
 \begin{equation}
\mbox{Find}\; \hat{p}\in {\cal A}\; \mbox{such that}\; J(\hat{p}) = 
\min_{{p}\in{\cal A}} J({p}).
   \label{minc}
 \end{equation}
Of course in real life $J$ does not vanish at the minimum because of
errors in the model and in the measurements.

Several questions should be adressed. Is the minimum unique~? Are
there local minima~? Is $\hat{p}$ very sensitive to uncertainties in
the measurements $z$~? In the next sections we shall present several
tools that are useful to answer these questions. They are based on
linear analysis (Sections \ref{sechess} and \ref{secedge}) and
nonlinear analysis (Sections \ref{secnl1} and \ref{secnl2}) of
stability.

But, before, let us consider the question of parameterization, always
crucial in parameter estimation.

\section{Parameterization}
\label{parameterization}

A very common choice for parameterization of the relative
permeabilities and capillary pressure is to use an analytical
representation:
\[ \begin{array}{ll}
k_w(S)=a_wS^{bw}, k_{nw}(S)=a_{nw}(1-S)^{bnw},\\
\displaystyle{\frac{dP_c}{dS}}=c_0+c_1S+c_2S^2,\;P_c(1.)=0.
\end{array}\]
The set of parameters to estimate is
$p=(a_w,b_w,a_{nw},b_{nw},c_0,c_1,c_2)$, a set of 7 parameters. With
such a choice it is usually not difficult to find a minimum to the
minimization problem (\ref{cout}),(\ref{minc}).

However there are many cases for which such a representation is not
suitable and the relative permeability and capillary curves do not
have such analytical shapes. To do without such a priori shapes a
possibility is to use a discrete representation of $kr_w, kr_{nw}, P_c$:
$p=(k_{rwj},j=1, k_{rnwj}, P_{cj}, j=1,\dots,n_s)$ where 
$k_{rwj},j=1, k_{rnwj}, P_{cj}$ are the values of $kr_w, kr_{nw}, P_c$
at a set of $n_s$ discretization points of the saturation interval
(0,1). Between these points the functions are linearly interpolated. 
Note that if one uses $n_s = 10$ equidistant saturation points, which
is reasonable to capture nonstandard shapes, then the number of
parameters is already 30.

As reported in \cite{chachajafliu90} for a standard displacement experiment, it
was possible to estimate the relative permeability curves using
cumulated production, pressure drop and saturation
measurements. However, without saturation measurement, the calculated
optimal parameters were depending on the initial guess of the
minimization algorithm. 
But when trying to estimate simultaneously relative permeability and
capillary pressure curves, the minimization algorithm would usually
get stuck in a local minimum with no practical interest.

\begin{figure}[htbp]
\begin{center}
\setlength{\unitlength}{0.666pt}
\begin{picture}(200,350)(0,0)
\includegraphics[height=8cm]{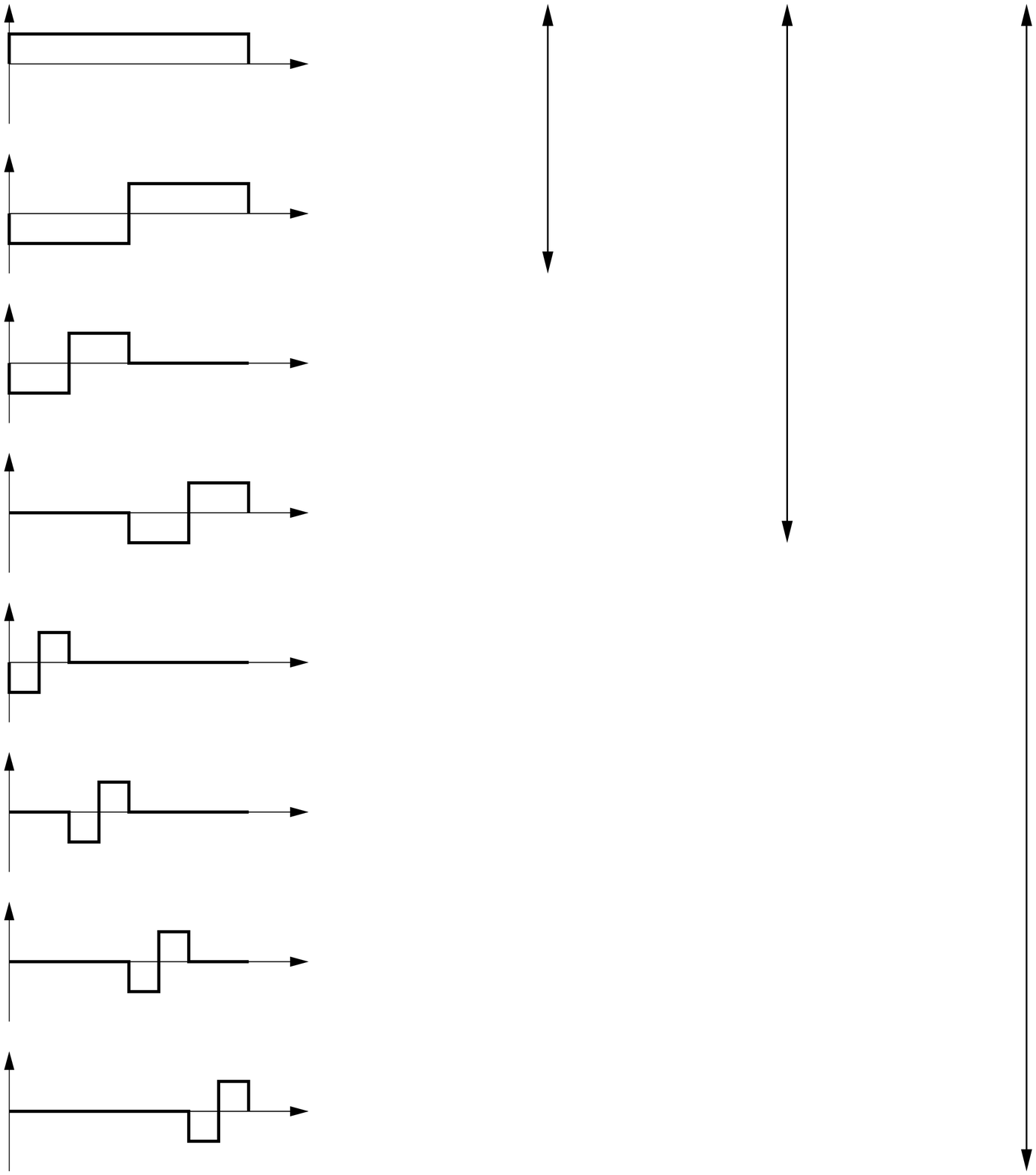}
\put(-338,320){\large $\Phi^0$}
\put(-307,317){{\scriptsize 0}}
\put(-231,314){{\scriptsize 1}}
\put(-205,320){Scale 0}
\put(-338,277){\large $\Phi^1_1$}
\put(-307,274){{\scriptsize 0}}
\put(-231,271){{\scriptsize 1}}
\put(-338,233){\large $\Phi^2_1$}
\put(-307,230){{\scriptsize 0}}
\put(-231,227){{\scriptsize 1}}
\put(-135,297){Scale 1}
\put(-338,189){\large $\Phi^2_2$}
\put(-307,186){{\scriptsize 0}}
\put(-231,183){{\scriptsize 1}}
\put(-65,252){Scale 2}
\put(-338,144){\large $\Phi^3_1$}
\put(-307,142){{\scriptsize 0}}
\put(-231,139){{\scriptsize 1}}
\put(-338,100){\large $\Phi^3_2$}
\put(-307,98){{\scriptsize 0}}
\put(-231,95){{\scriptsize 1}}
\put(-338,56){\large $\Phi^3_4$}
\put(-307,54){{\scriptsize 0}}
\put(-231,51){{\scriptsize 1}}
\put(-338,12){\large $\Phi^3_4$}
\put(-307,10){{\scriptsize 0}}
\put(-231,7){{\scriptsize 1}}
\put(10,160){Scale 3}
\end{picture}
\end{center}
\caption{The Haar basis}
\label{haar}
\end{figure}

This is the reason for the introduction of multiscale
parameterization. The main idea is to adapt the parameterization as
the minimization advances, starting with few parameters in order to
estimate the main features of the functions, and increasing slowly their
number to estimate their refined features. A simple example of a
multiscale basis to expand a function is the Haar basis as represented
in Fig. \ref{haar}.

In practice one endpoint of each of the relative permeability and
capillary curves is known, so the continuous piecewise linear
representation of the functions are uniquely defined by their
piecewise constant derivatives which will be parameterized with the
Haar basis :
\[
\displaystyle{\frac{dk_{rw}(S)}{dS}=(c_rw)^0\Phi^0(S)+\sum^{n_s-1}_{i=1}\sum^{2^{i-1}}_{j=1}
(c_rw)^i_j\Phi^i_j(S)}, \]
and similarly for $\displaystyle{\frac{dk_{rnw}(S)}{dS},
\frac{dP_c(S)}{dS}}$. Note that with such a parameterization enforcing
conditions (\ref{properties}) is simple. 

\begin{figure}[htbp]
\includegraphics[height=4cm]{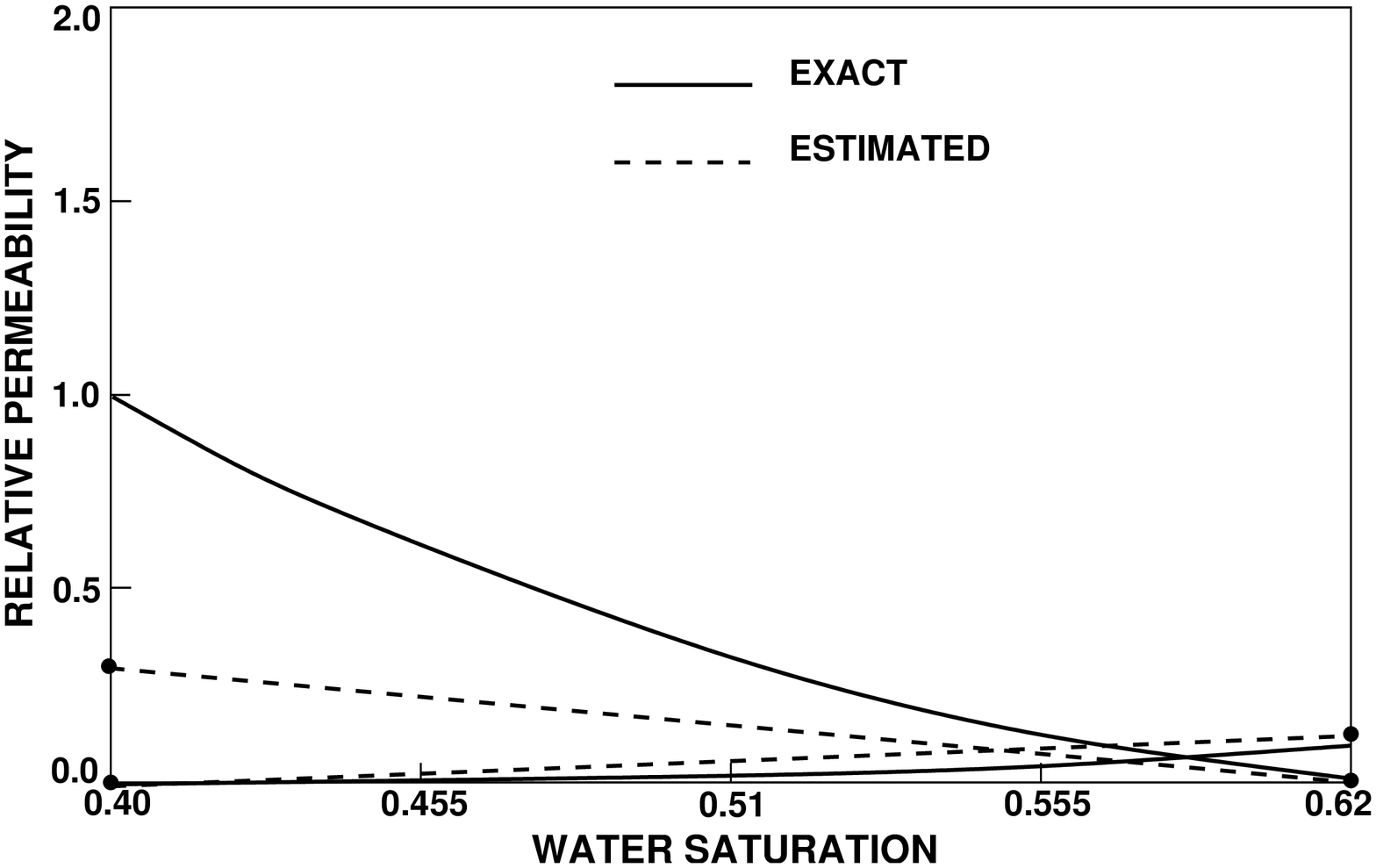} \raisebox{2cm}{{\bf Scale 0}}
\includegraphics[height=4cm]{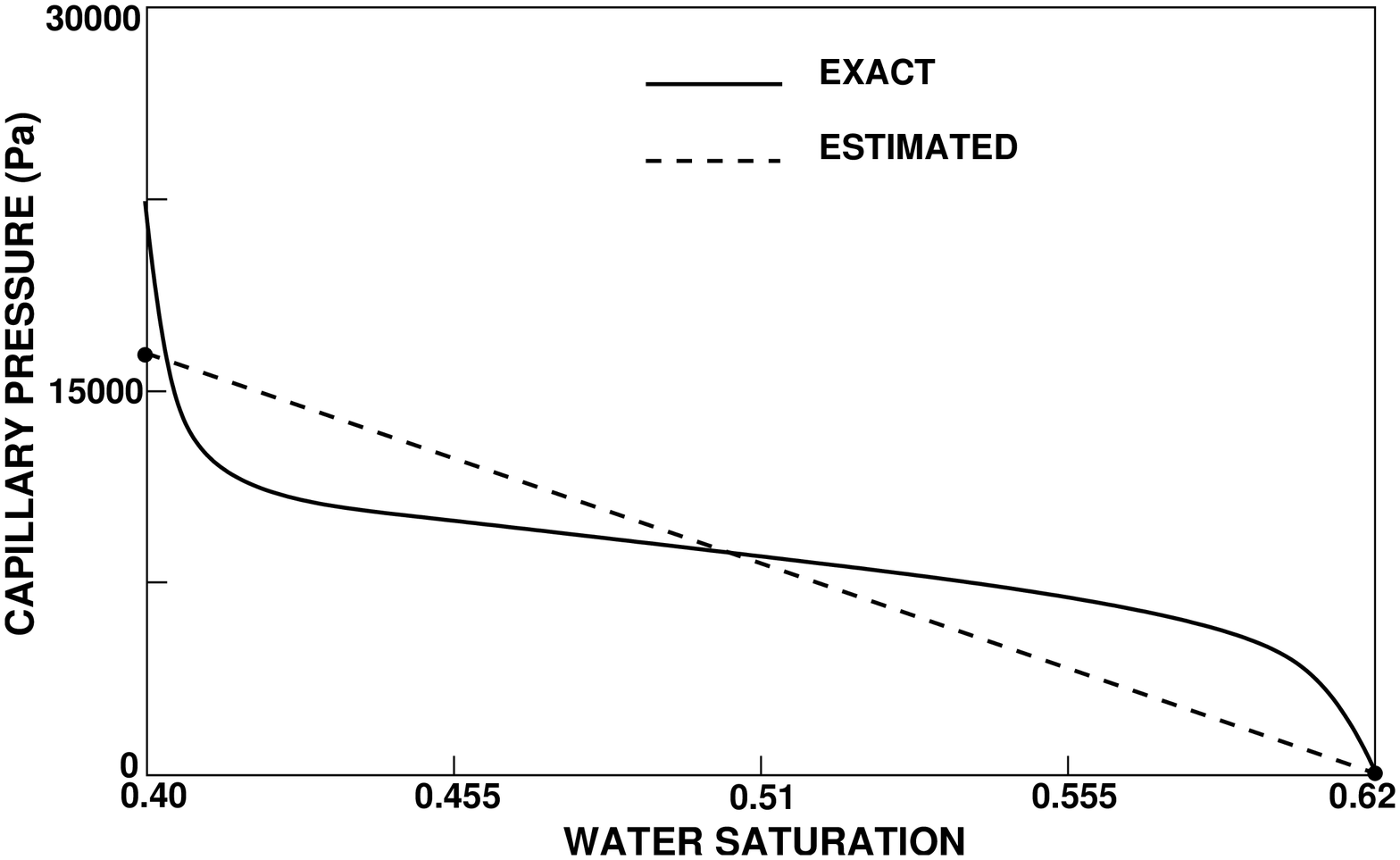} \\
\includegraphics[height=4cm]{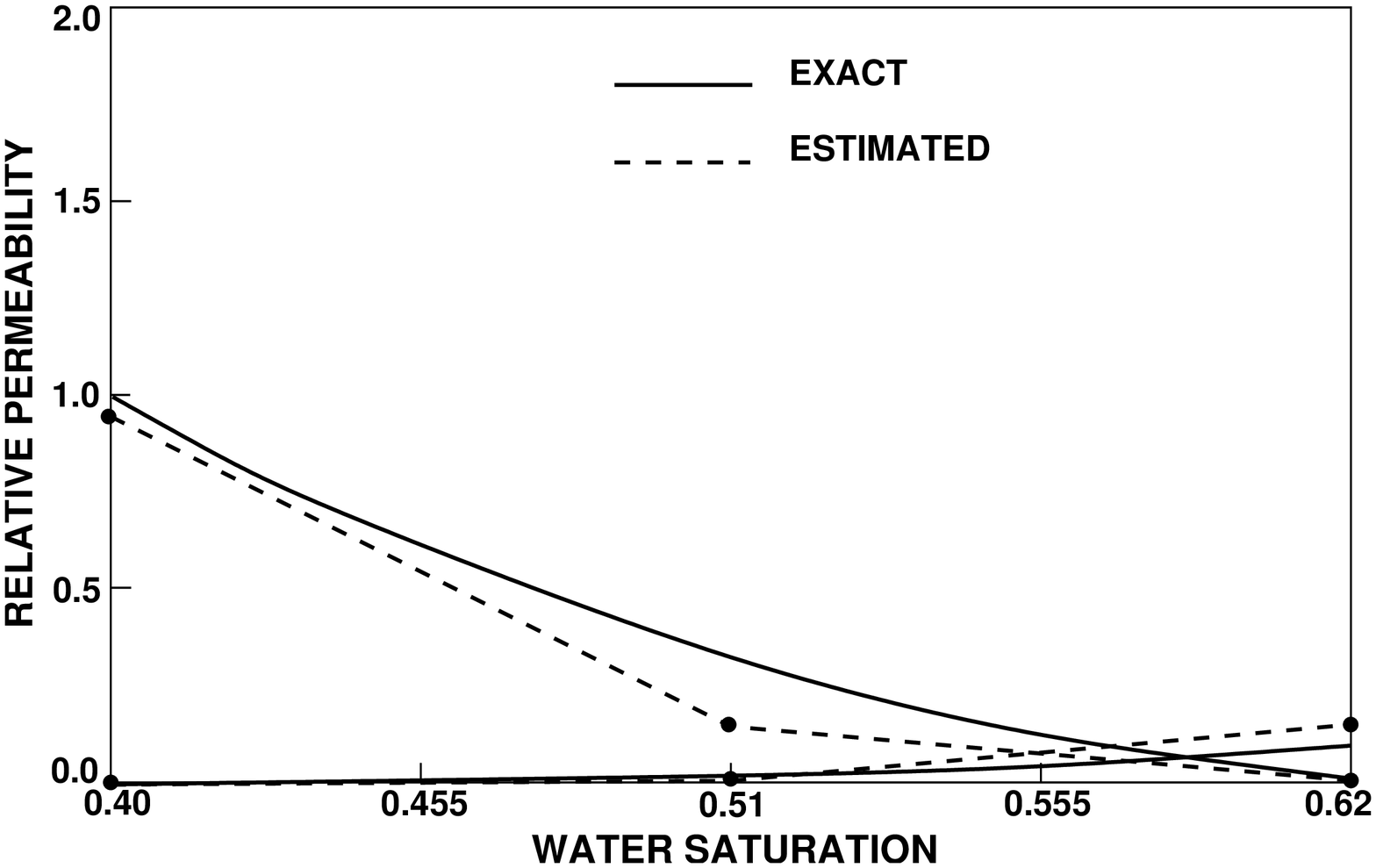} \raisebox{2cm}{{\bf Scale 1}}
\includegraphics[height=4cm]{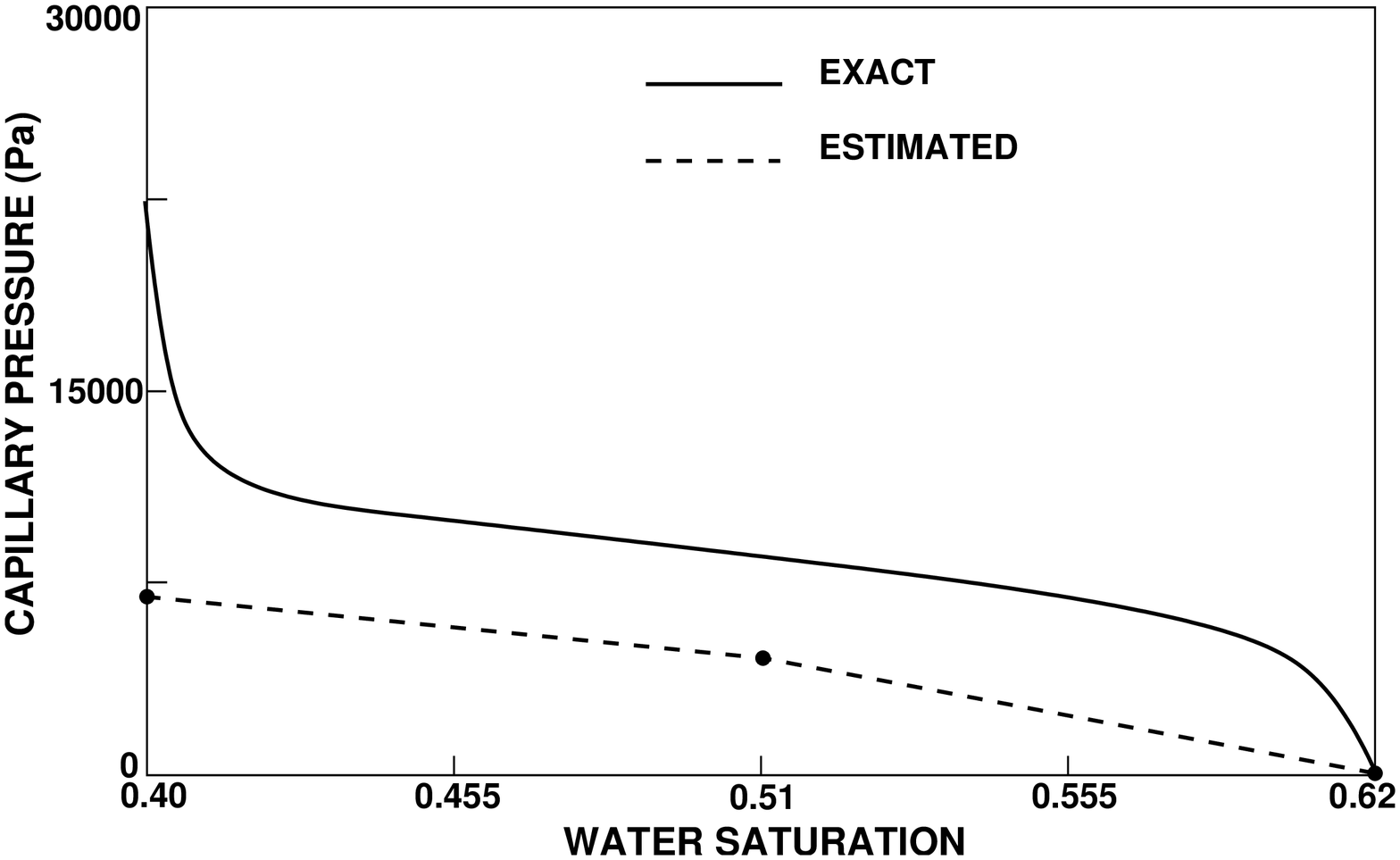} \\
\includegraphics[height=4cm]{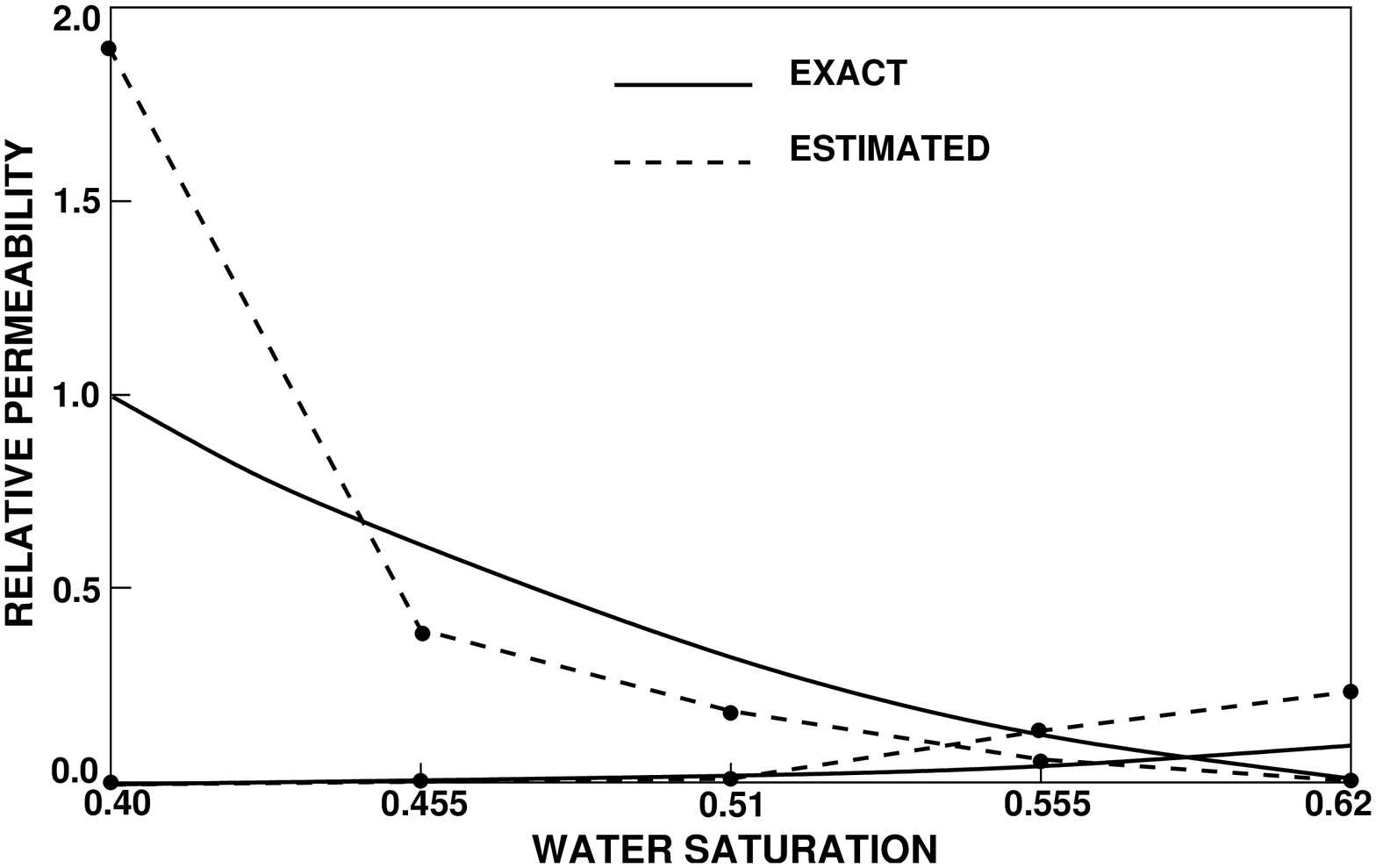} \raisebox{2cm}{{\bf Scale 2}}
\includegraphics[height=4cm]{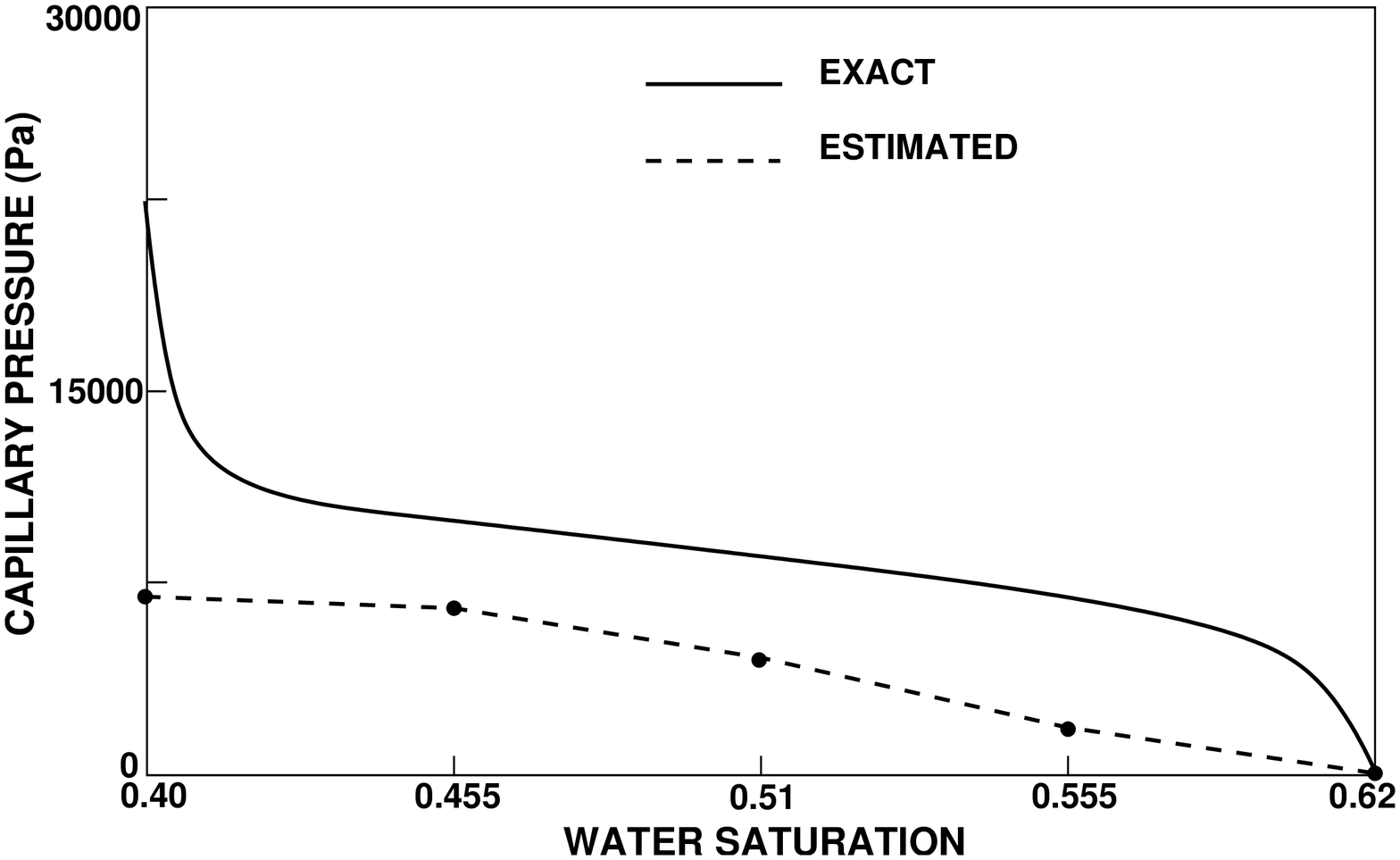} \\
\includegraphics[height=4cm]{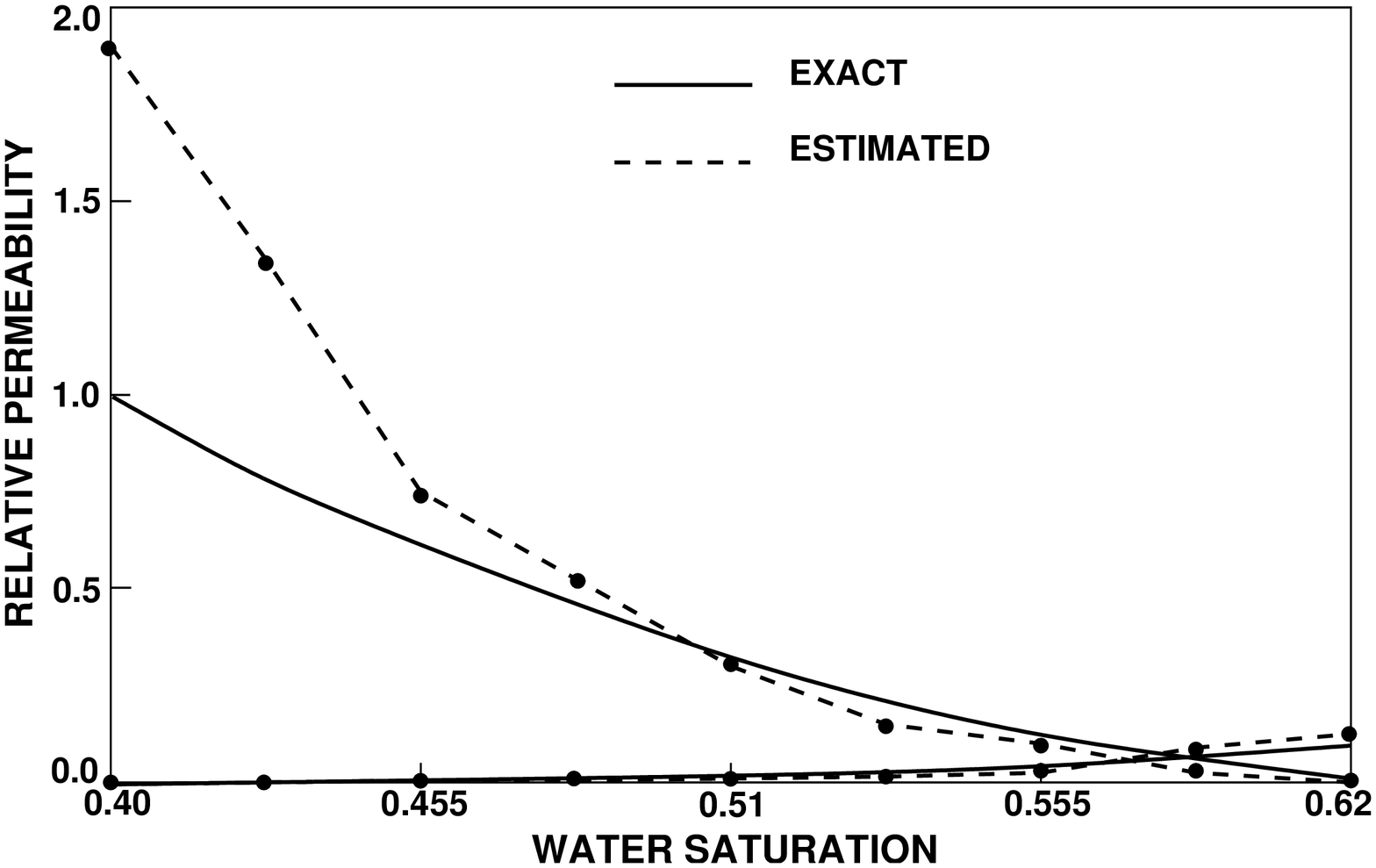} \raisebox{2cm}{{\bf Scale 3}}
\includegraphics[height=4cm]{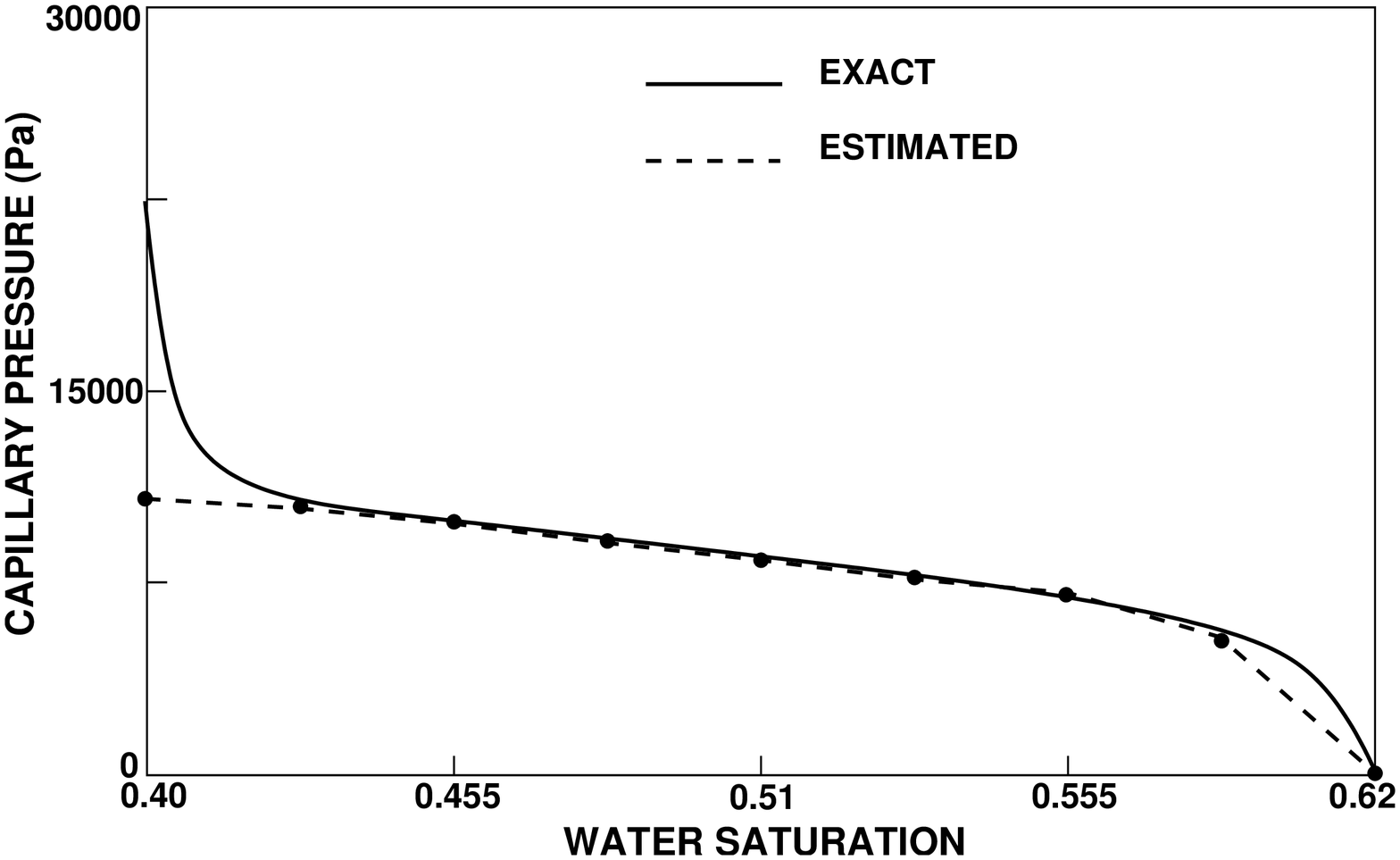} \\
\includegraphics[height=4cm]{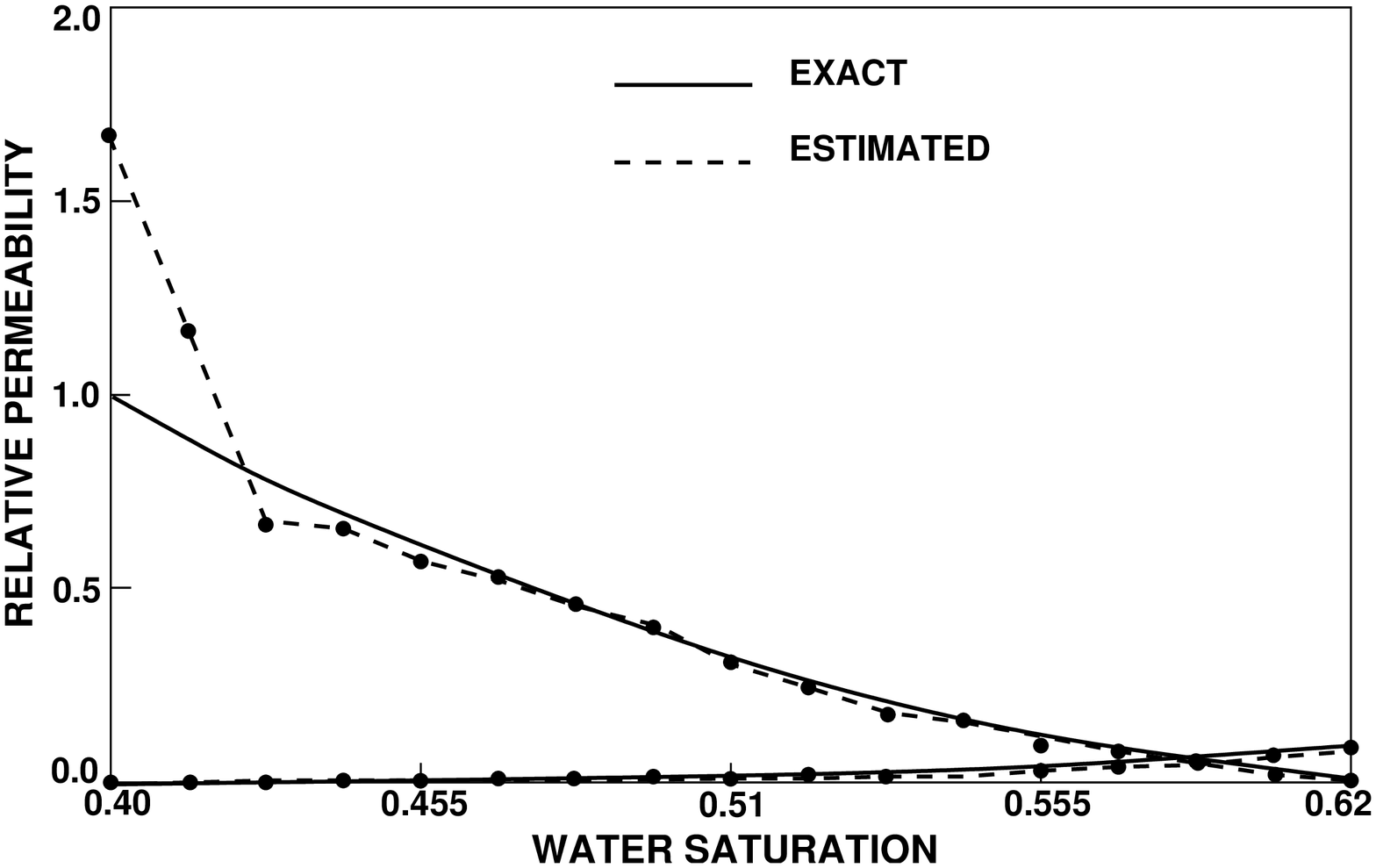} \raisebox{2cm}{{\bf Scale 4}}
\includegraphics[height=4cm]{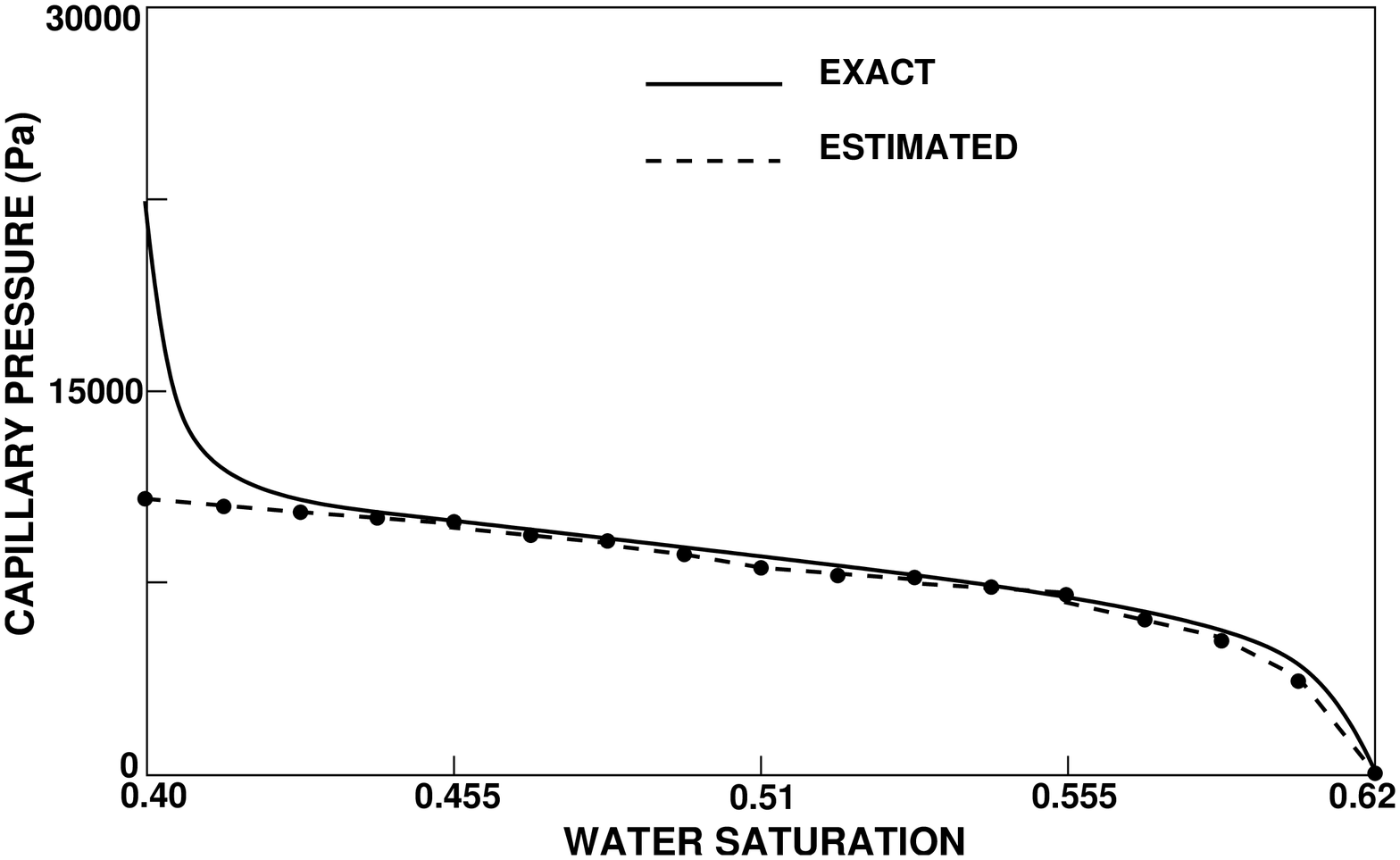} \\
\caption{Relative permeability and capillary pressure curves estimated
with multiscale parameterization}
\label{multiresults}
\end{figure}

Then a multiscale optimization proceeds as follows:
\begin{enumerate}
\item Minimize at scale $i=0$.
\item Augment $i$ by $1$.
\item Minimize at scale $i$ using as starting point functions obtained by
linearly interpolating that estimated with scale $i-1$.
\item If $i\geq i_{\max}$ or if scales $i-1$ and $i$ give the same estimated
function, stop.\\ If not go to 2.
\end{enumerate}
With such a procedure, the simultaneous estimation of relative
permeability and capilary curves was successful while it failed when
using a standard discrete parameterization \cite{chachajafliu90}. Fig.~\ref{multiresults}
shows the progression of the multiscale parameterization for relative
permeabilities and capillary pressure.

An important advantage of multiscale optimization is that it is
adaptive. It does not require to fix a priori the number of parameters
which increases as the number of steps of multiscale optimization
increases until the difference between two scales is not anymore
significant. Therefore the multiscale optimization procedure
determines itself an almost optimal number of parameters necessary to
interpret the available data.
\section{Stability analysis based on the Hessian}
\label{sechess}
In this section we show how to obtain some information concerning the
inverse problem using the Hessian of $J$. We linearize $\varphi$
around a given parameter $p_0$. Let us perturbate $p_0$ into $p$ with
a small perturbation $\delta {p}={p}-{p}_0$.
Taylor's expansion gives
\[  \varphi(p) \approx \varphi(p_0)+{\varphi}'(p_0)\delta p. \]
If we introduce $\delta z =z -\varphi(p_0)$, the observation error function can actually be written as a quadratic function of
$\delta p$ :
\begin{equation} \label{defjl}
J(\delta p)=\parallel{\varphi}'(p_0)\delta p-\delta z
{\parallel}^2_W.
\end{equation}

Then the inverse problem is to minimize this function and
the minimum $\delta {\hat p}$ satisfies   
\begin{equation}
 H \delta {\hat p}=  {{\varphi}'(p_0)}^t W \delta z 
  \label{eqhessien}
\end{equation}
with $H = {{\varphi}'(p_0)}^t W {\varphi}'(p_0)$ the Hessian. We
note that the matrix ${\varphi}'(p_0)$, which is the jacobian matrix
of $\varphi$ at $p_0$, can be viewed as the sensitivity matrix for
$\varphi$. 

Therefore solving the linearized inverse problem reduces to solving
the linear system (\ref{eqhessien}) and this explains the importance of
studying the Hessian $H$.

To illustrate this we shall consider three problems of parameter
estimation for which, from realistic data obtained in centrifugal experiments \cite{chaforzhachalen92}, we generated the measured
observations with our simulation code. These problems are :
\begin{enumerate}
\item estimating relative permeabilities while measuring saturation profiles,
\item estimating relative permeabilities while measuring productions,
\item estimating relative permeabilities and capillary pressure while measuring
saturation profiles.
\end{enumerate}
The aim of this analysis is to study the importance of the choice
of measurements for estimating $kr_w, kr_{nw}, p_c$.

\begin{figure}[htp]
\begin{center}
\begin{minipage}[t]{6cm}
\begin{center}
\epsfysize=5cm \leavevmode 
\hspace*{-1cm}\epsffile{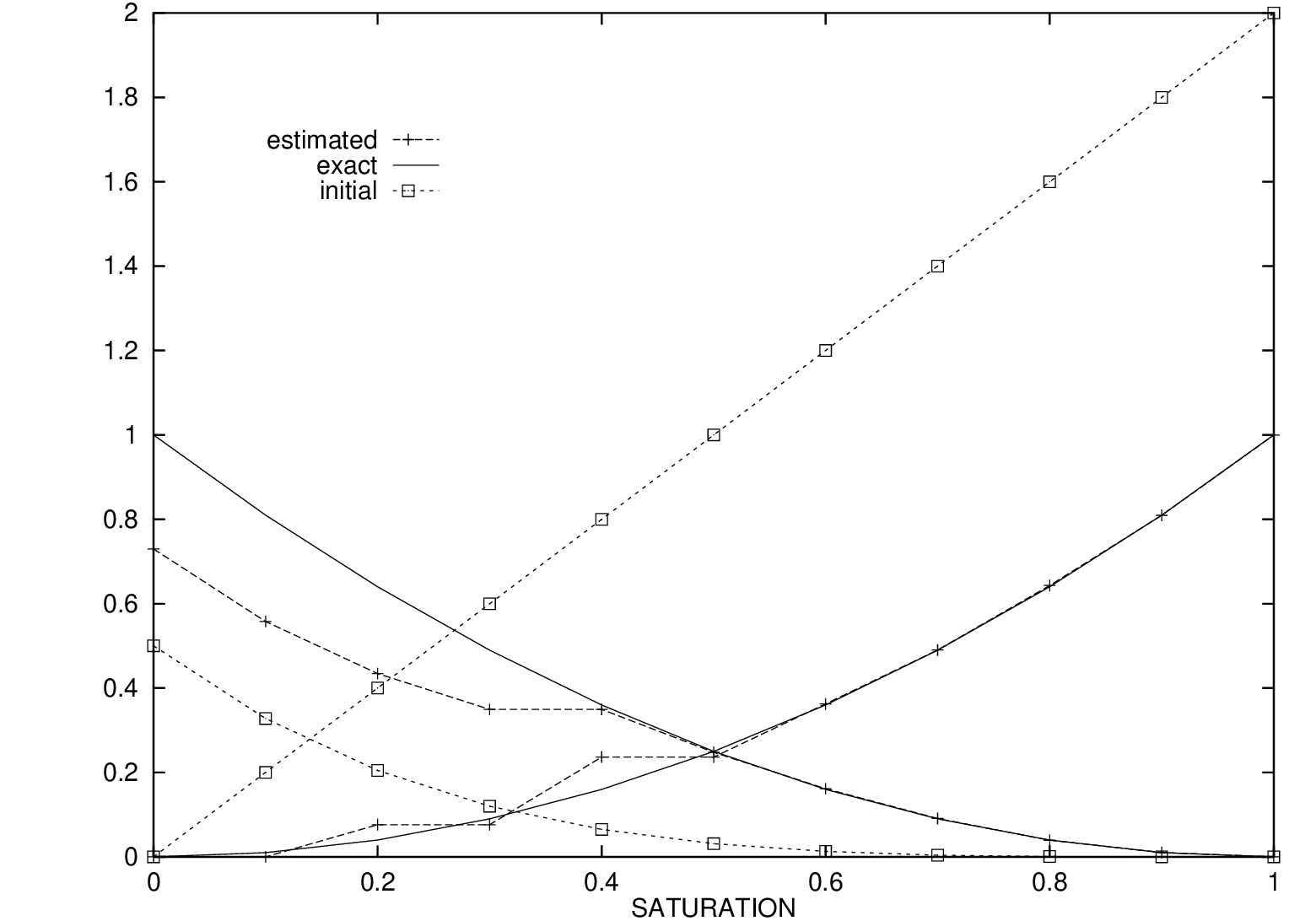}\\
Relative permeabilities 
\end{center}
\end{minipage}
\begin{minipage}[t]{6cm}
\begin{center}
\epsfysize=5cm \leavevmode 
\epsffile{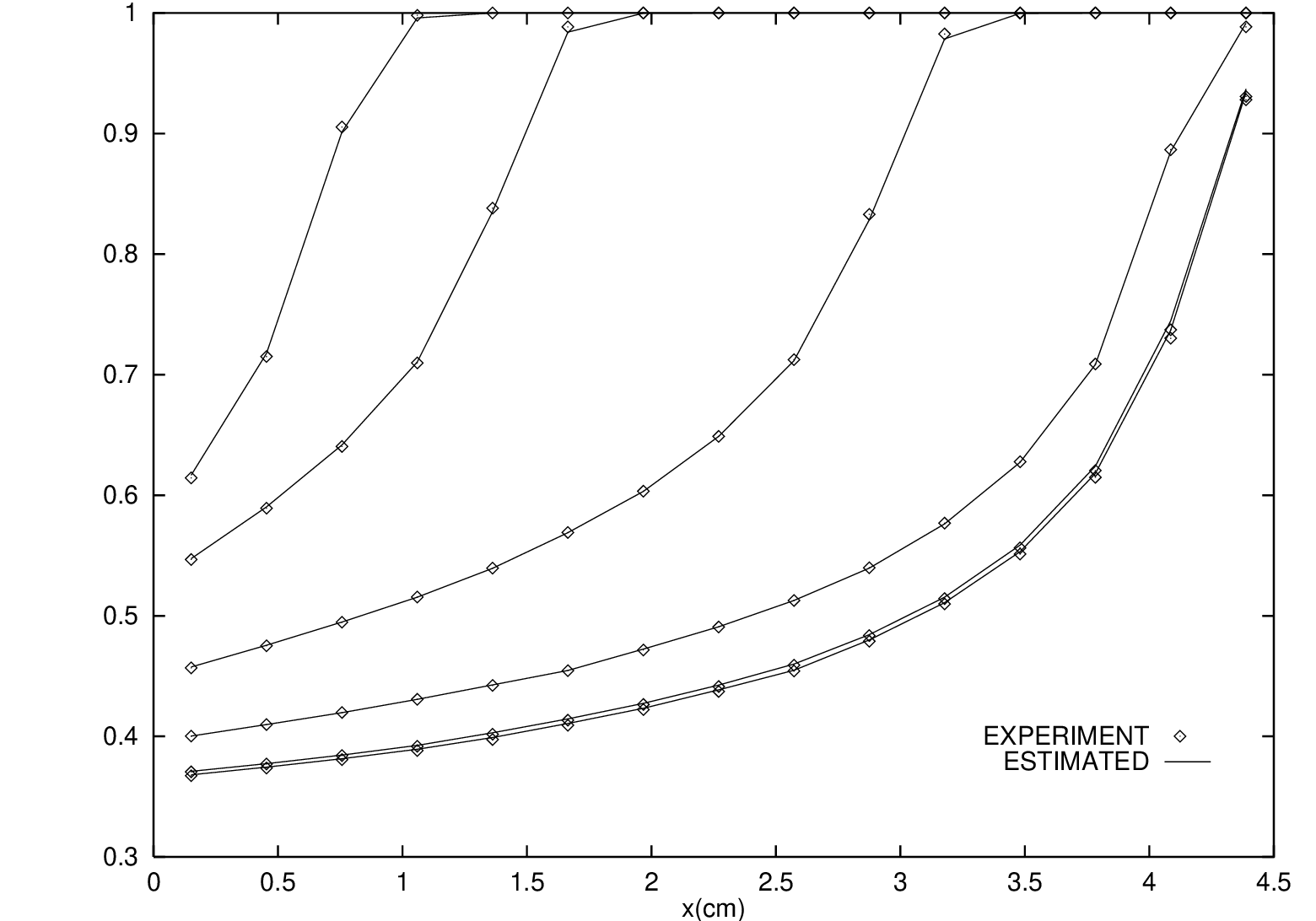}\\
\hspace*{1cm} Saturation profiles
\end{center}
\end{minipage}\\
\epsfysize=6cm \leavevmode 
\epsffile{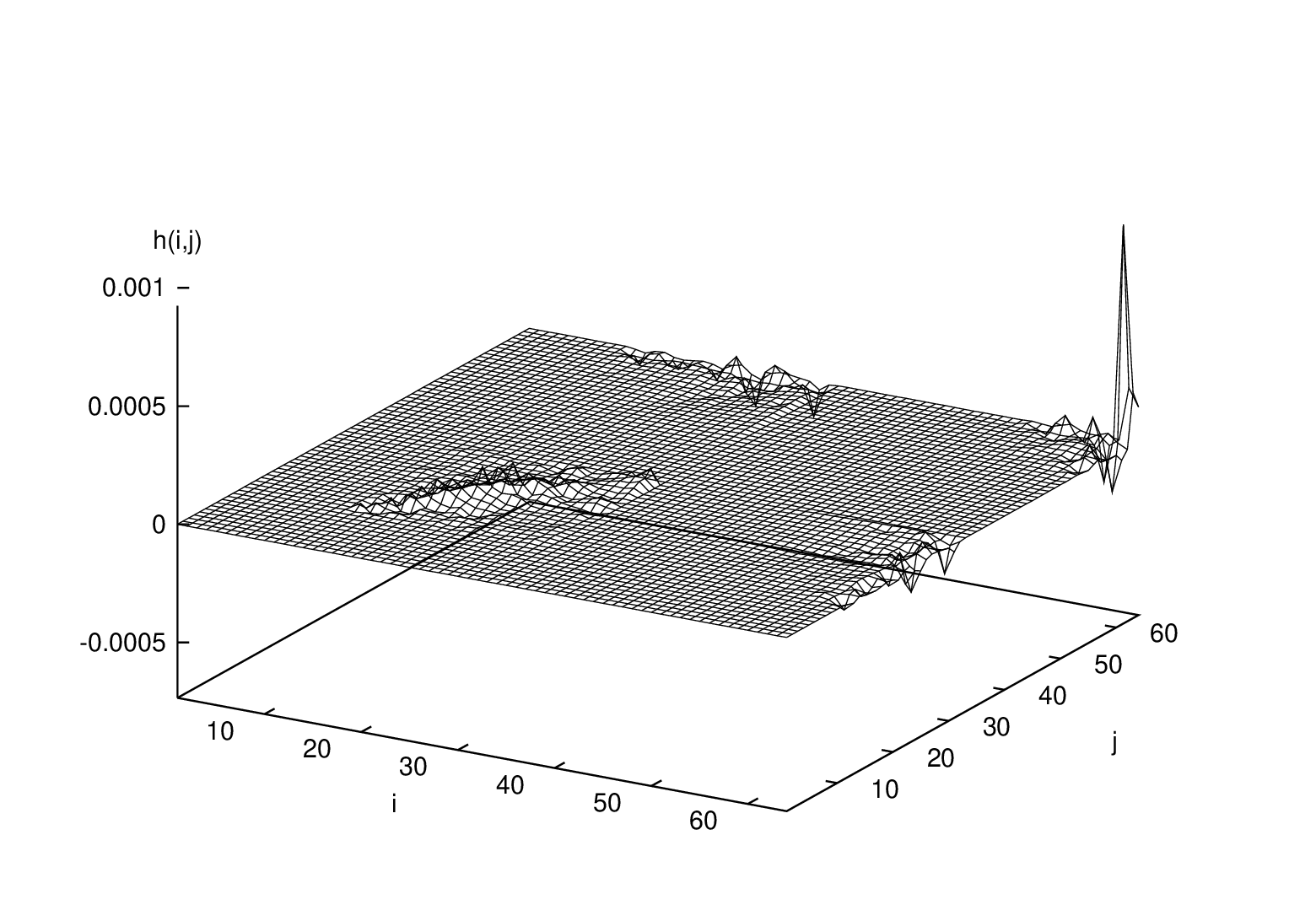}\\
The Hessian $H$\\[0.5cm]
\begin{minipage}[t]{6cm}
\begin{center}
\epsfysize=5cm \leavevmode 
\hspace*{-1cm}\epsffile{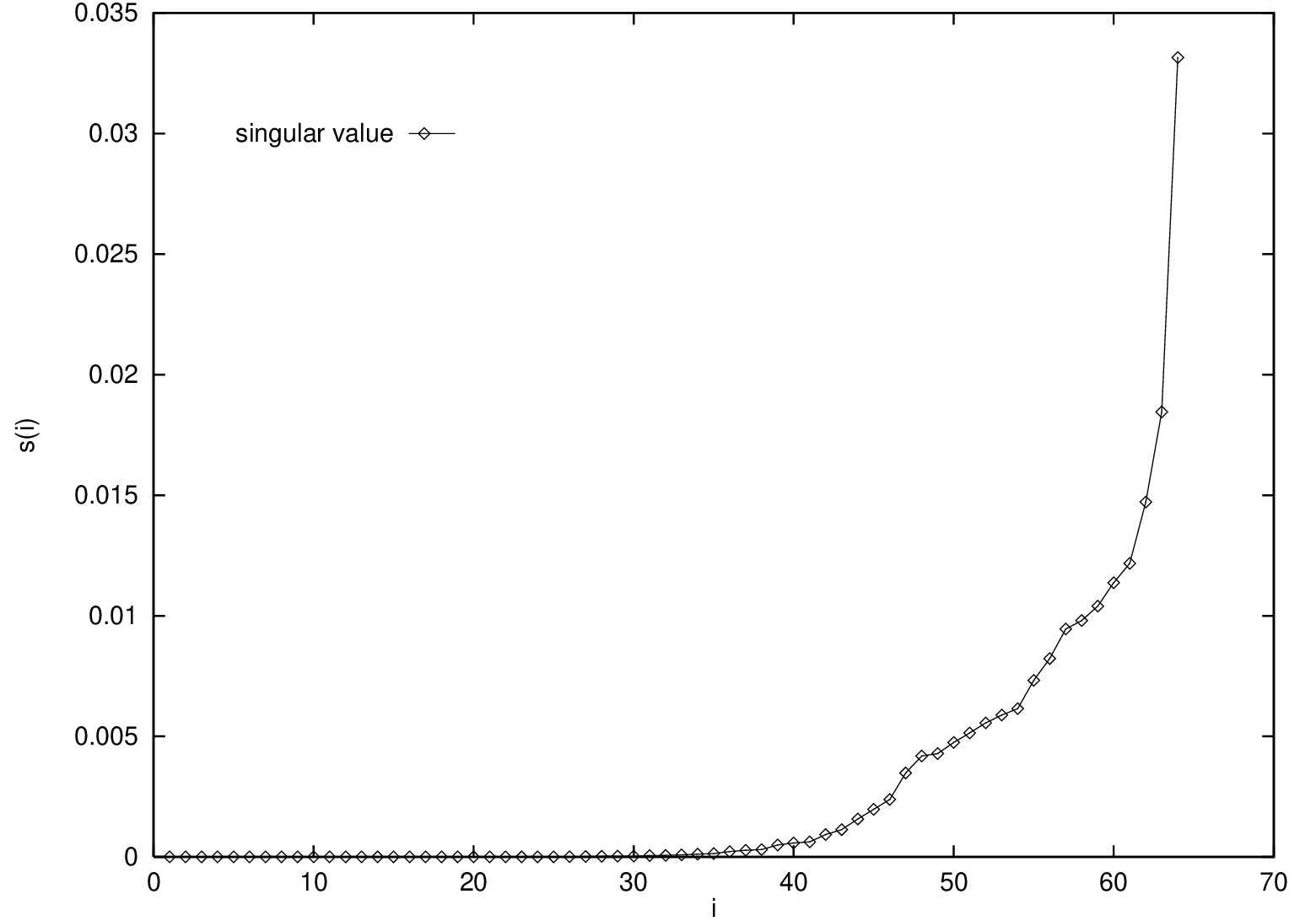}\\
Singular Values 
\end{center}
\end{minipage}
\begin{minipage}[t]{6cm}
\begin{center}
\epsfysize=5cm \leavevmode 
\epsffile{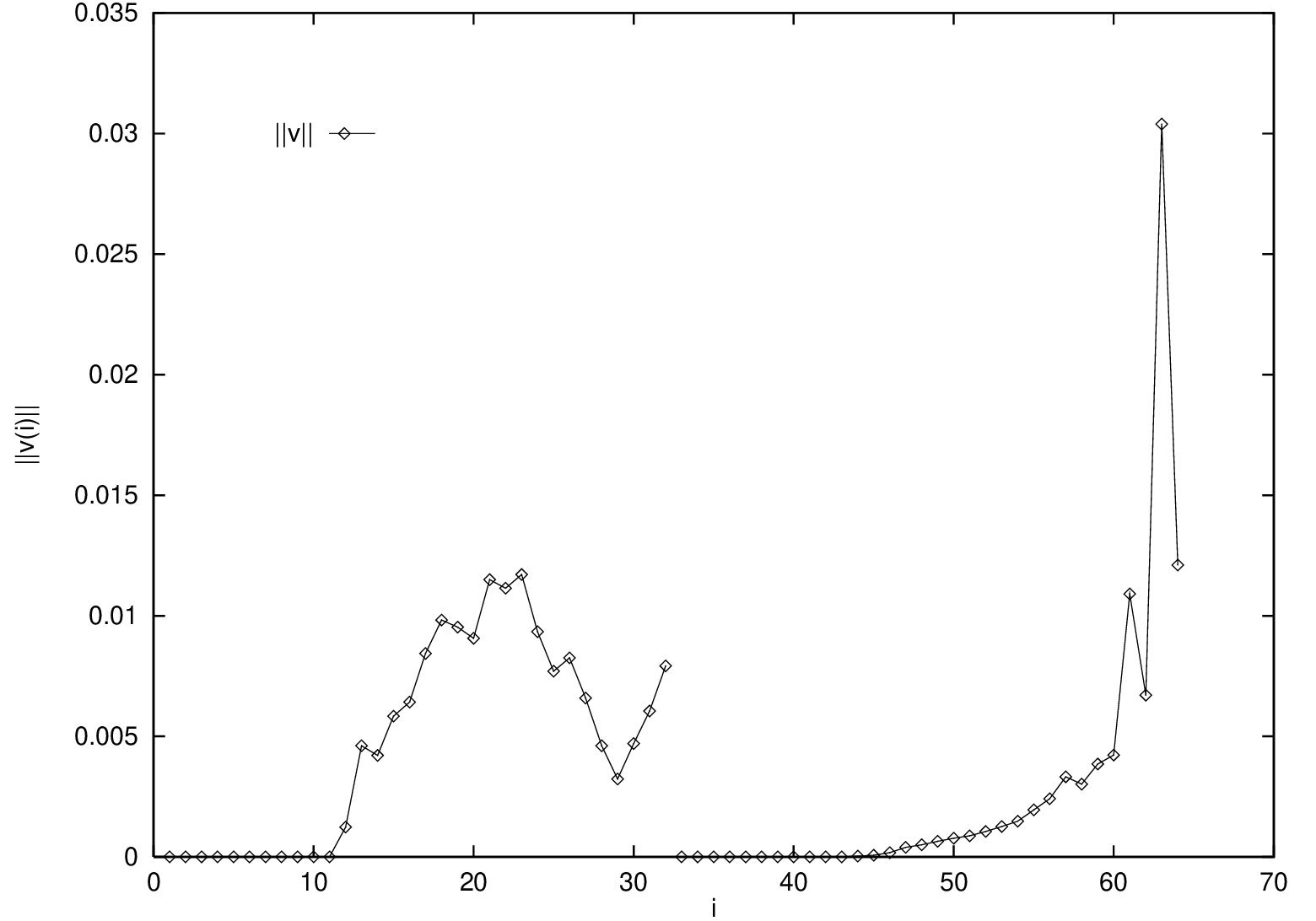}\\
Sensitivity of saturation profiles
\end{center}
\end{minipage}
\end{center}
\caption{Problem 1 : Estimating relative permeabilities while
  measuring saturation profiles}
\label{experiment1}
\end{figure}

\begin{figure}[htp]
\begin{center}
\begin{minipage}[t]{6cm}
\begin{center}
\epsfysize=5cm \leavevmode 
\hspace*{-1cm}\epsffile{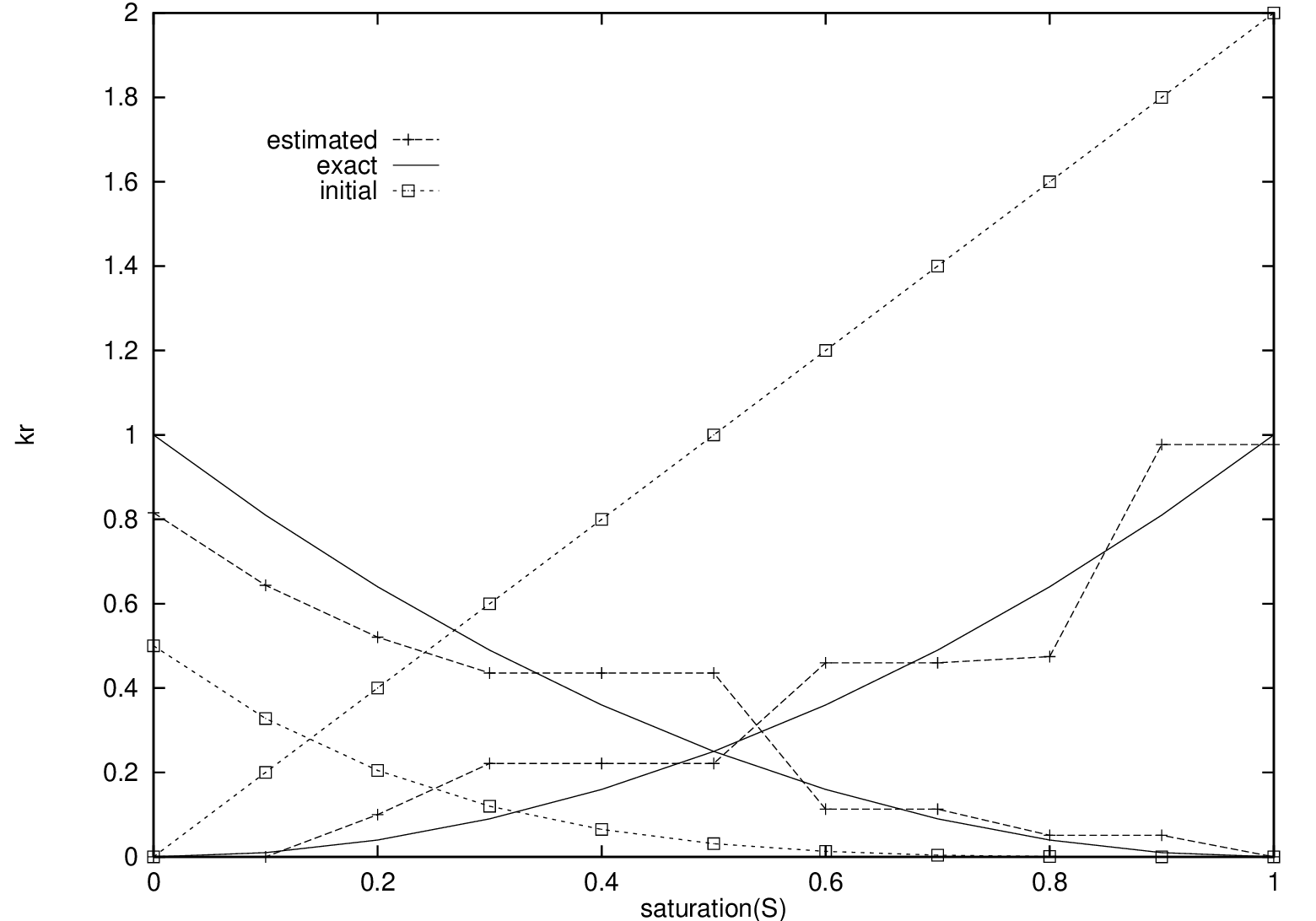}\\
Relative permeabilities 
\end{center}
\end{minipage}
\begin{minipage}[t]{6cm}
\begin{center}
\epsfysize=5cm \leavevmode 
\epsffile{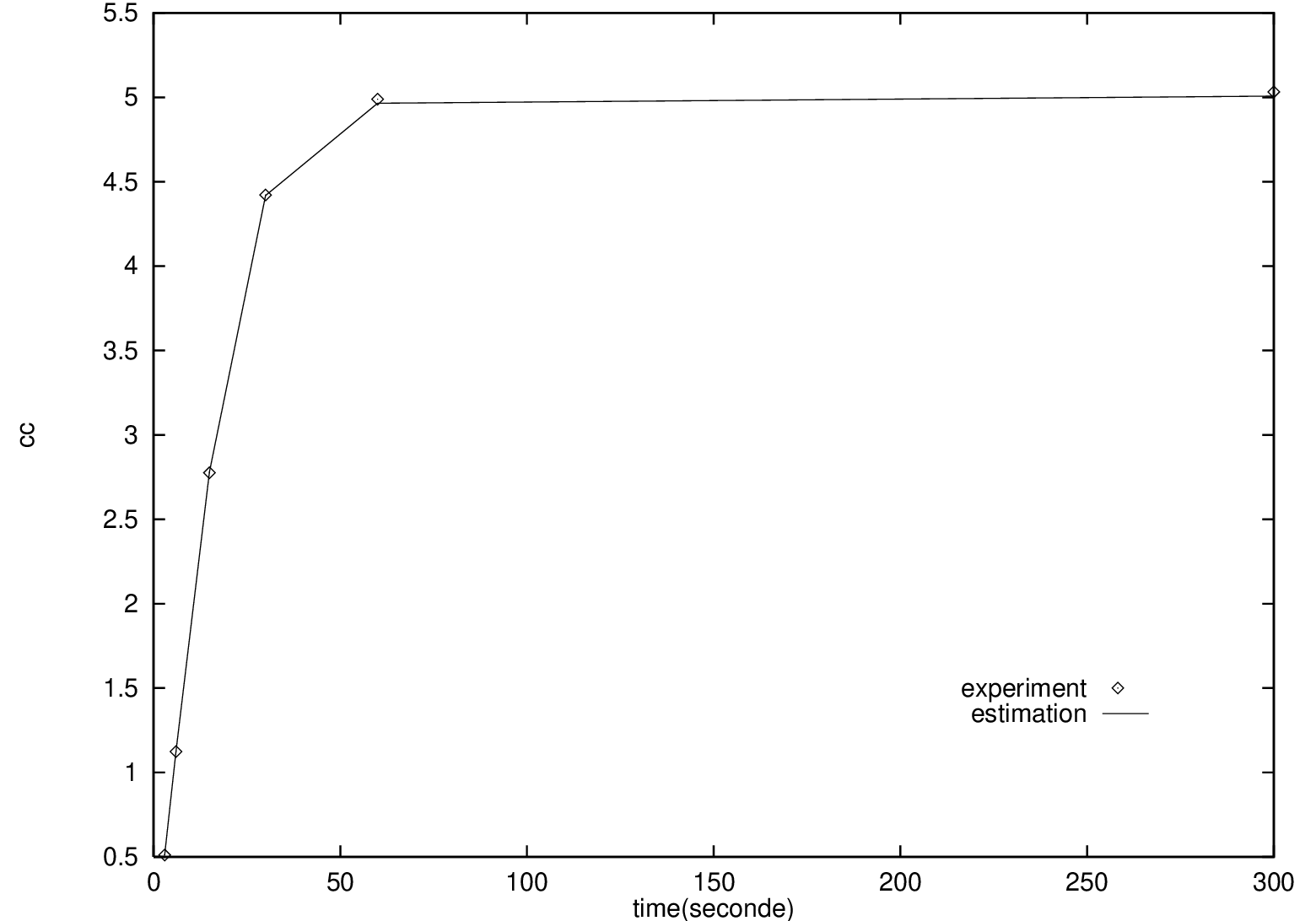}\\
\hspace*{1cm} Productions
\end{center}
\end{minipage}
\epsfysize=6cm \leavevmode 
\epsffile{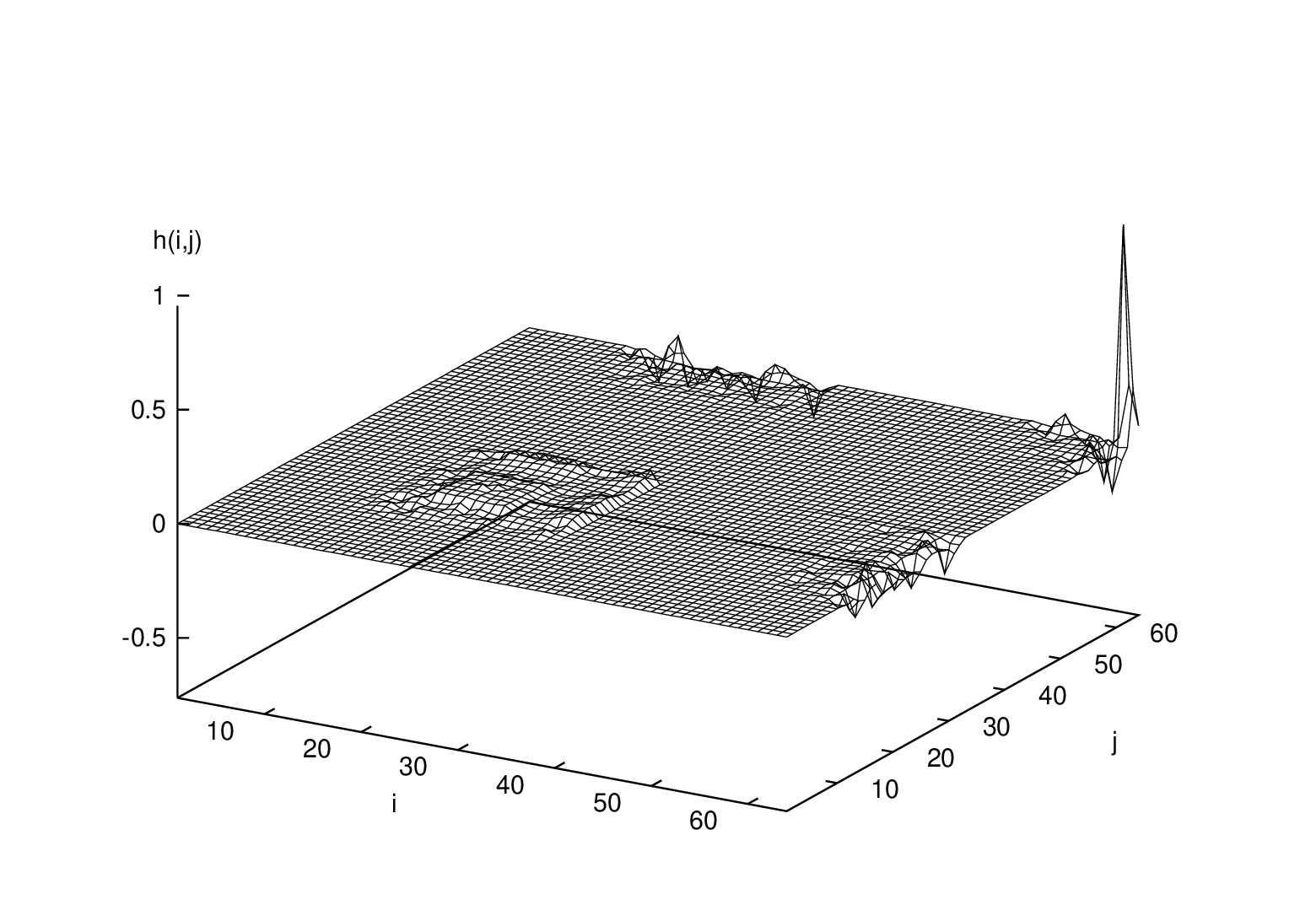}\\
The Hessian $H$\\[0.5cm]
\begin{minipage}[t]{6cm}
\begin{center}
\epsfysize=5cm \leavevmode 
\hspace*{-1cm}\epsffile{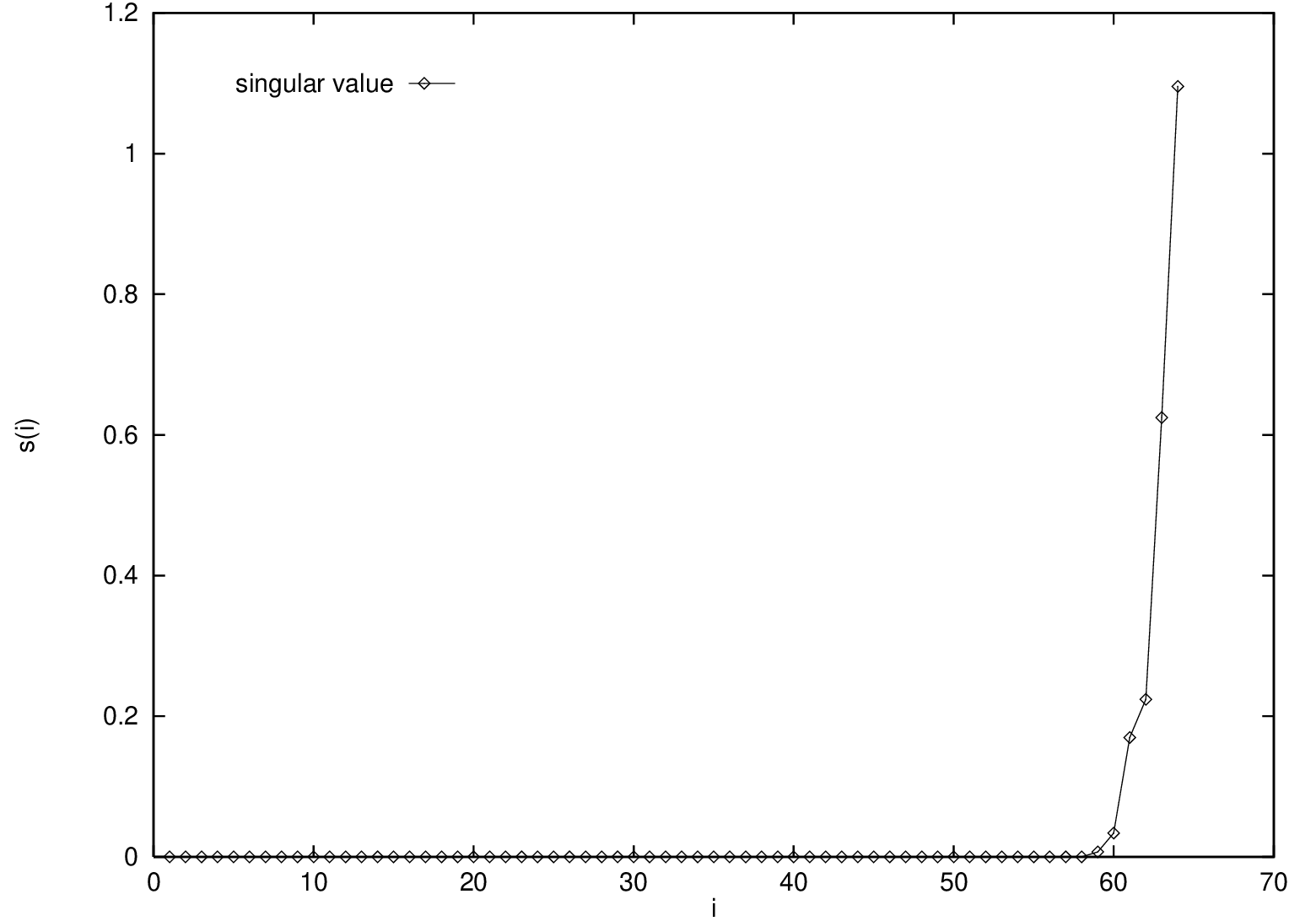}\\
Singular Values 
\end{center}
\end{minipage}
\begin{minipage}[t]{6cm}
\begin{center}
\epsfysize=5cm \leavevmode 
\epsffile{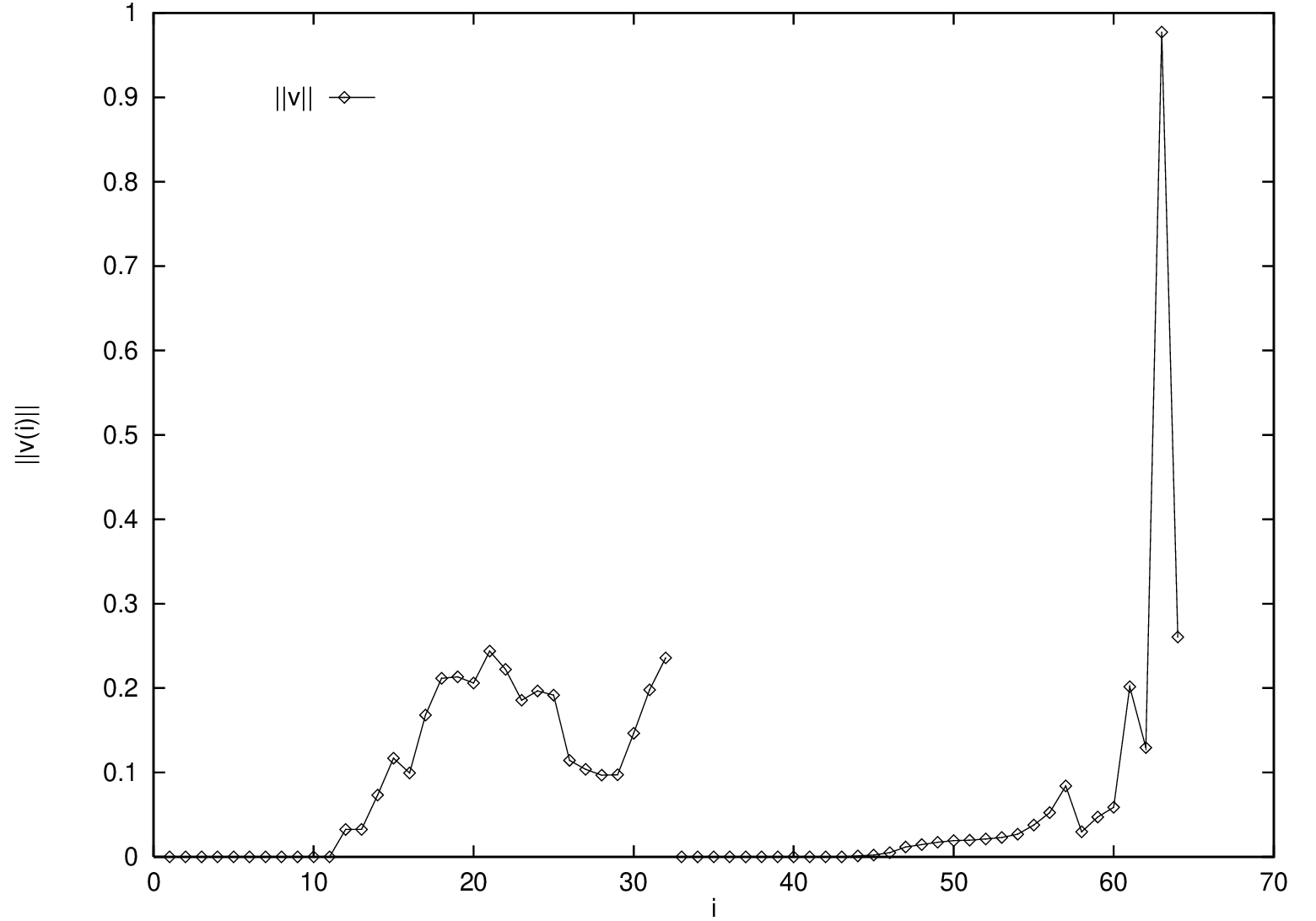}\\
Sensitivity of saturation profiles
\end{center}
\end{minipage}
\end{center}
\caption{Problem 2 : Estimating relative permeabilities while
  measuring productions}
\label{experiment2}
\end{figure}

\begin{figure}[htp]
\begin{center}
\begin{minipage}[t]{6cm}
\begin{center}
\epsfysize=5cm \leavevmode 
\hspace*{-1cm}\epsffile{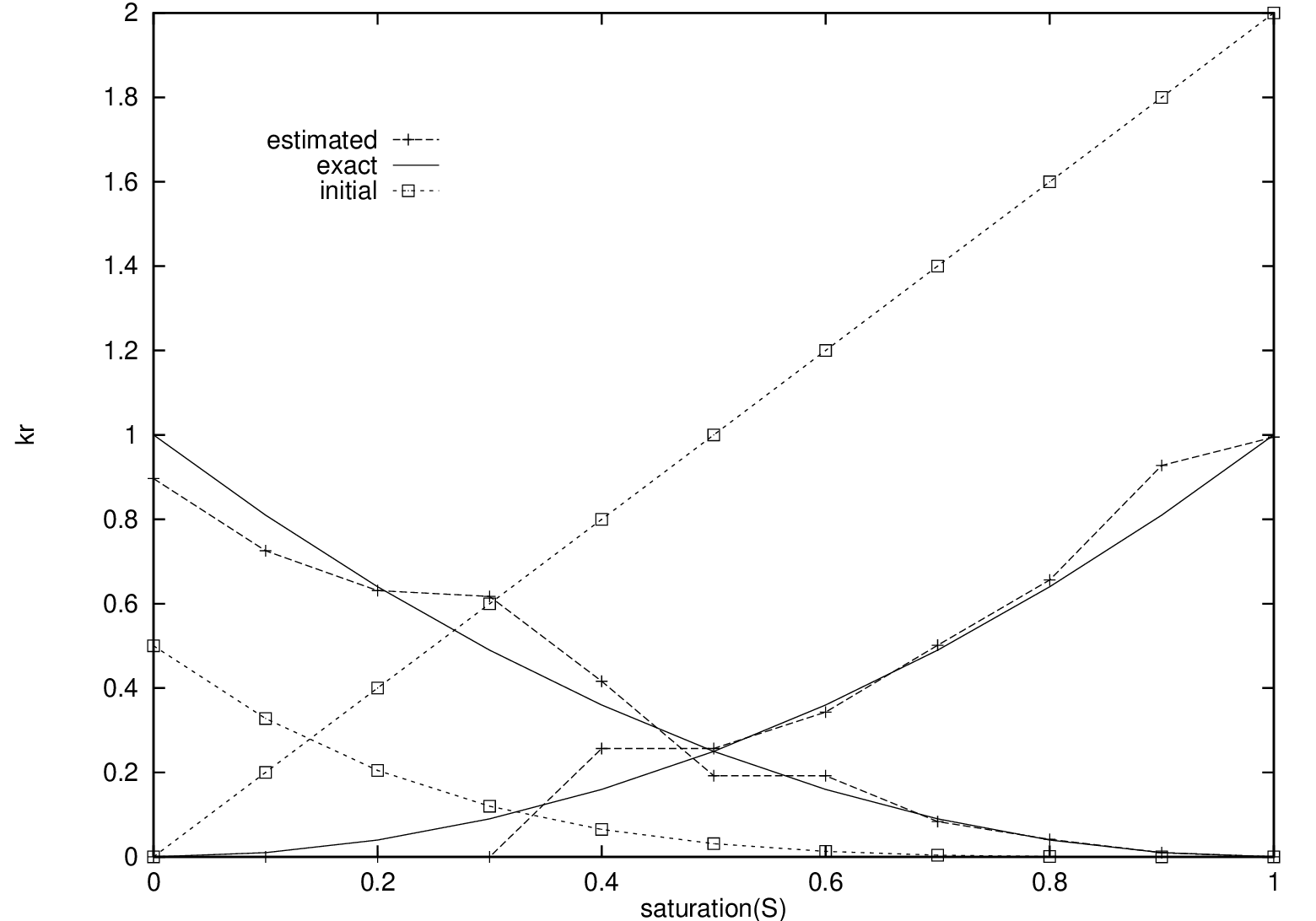}\\
Relative permeabilities 
\end{center}
\end{minipage} \hspace*{1cm}
\begin{minipage}[t]{6cm}
\begin{center}
\epsfysize=5cm \leavevmode 
\hspace*{-1cm}\epsffile{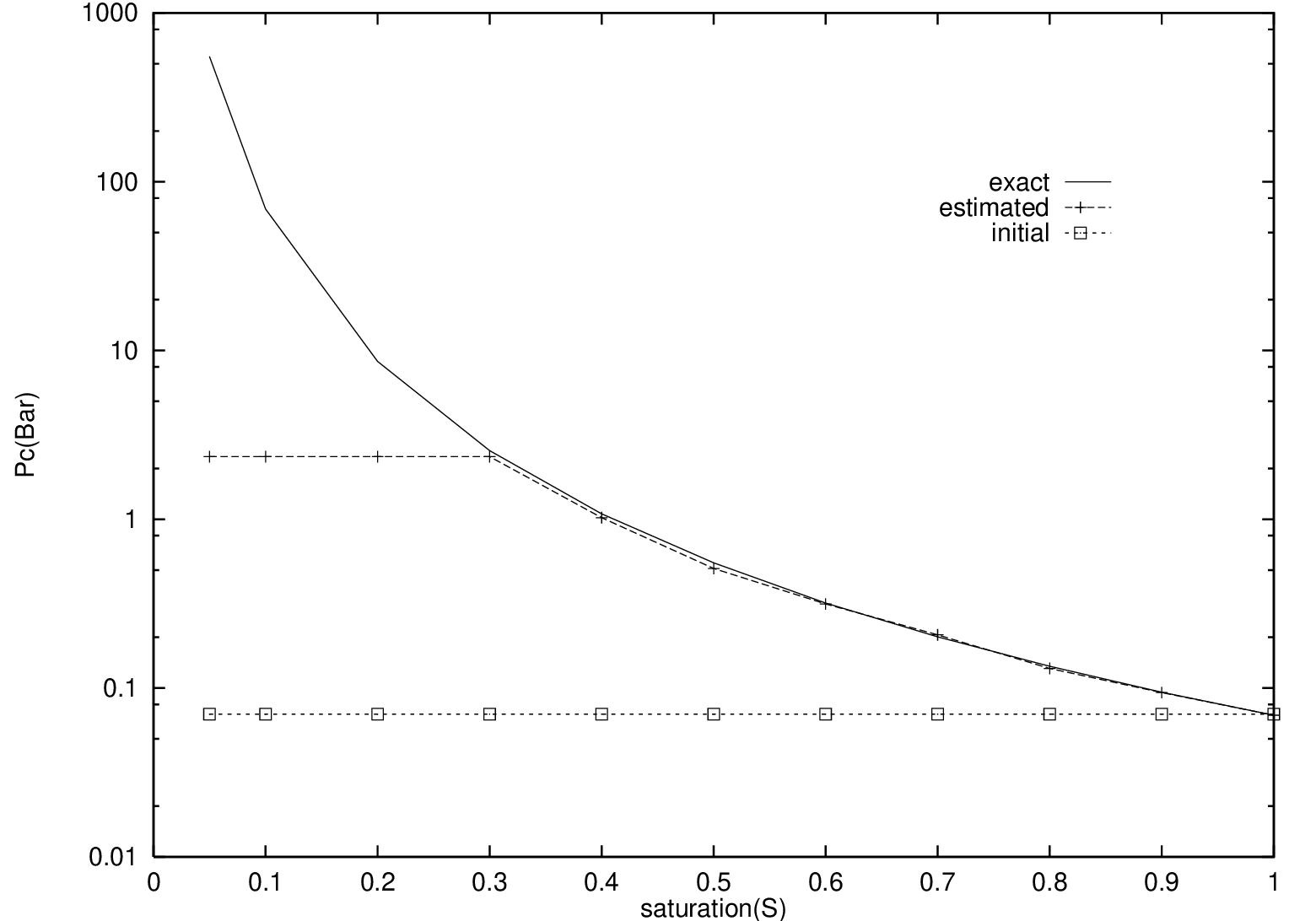}\\
Capillary pressure
\end{center}
\end{minipage}\\[0.5cm]
\epsfysize=5cm \leavevmode 
\epsffile{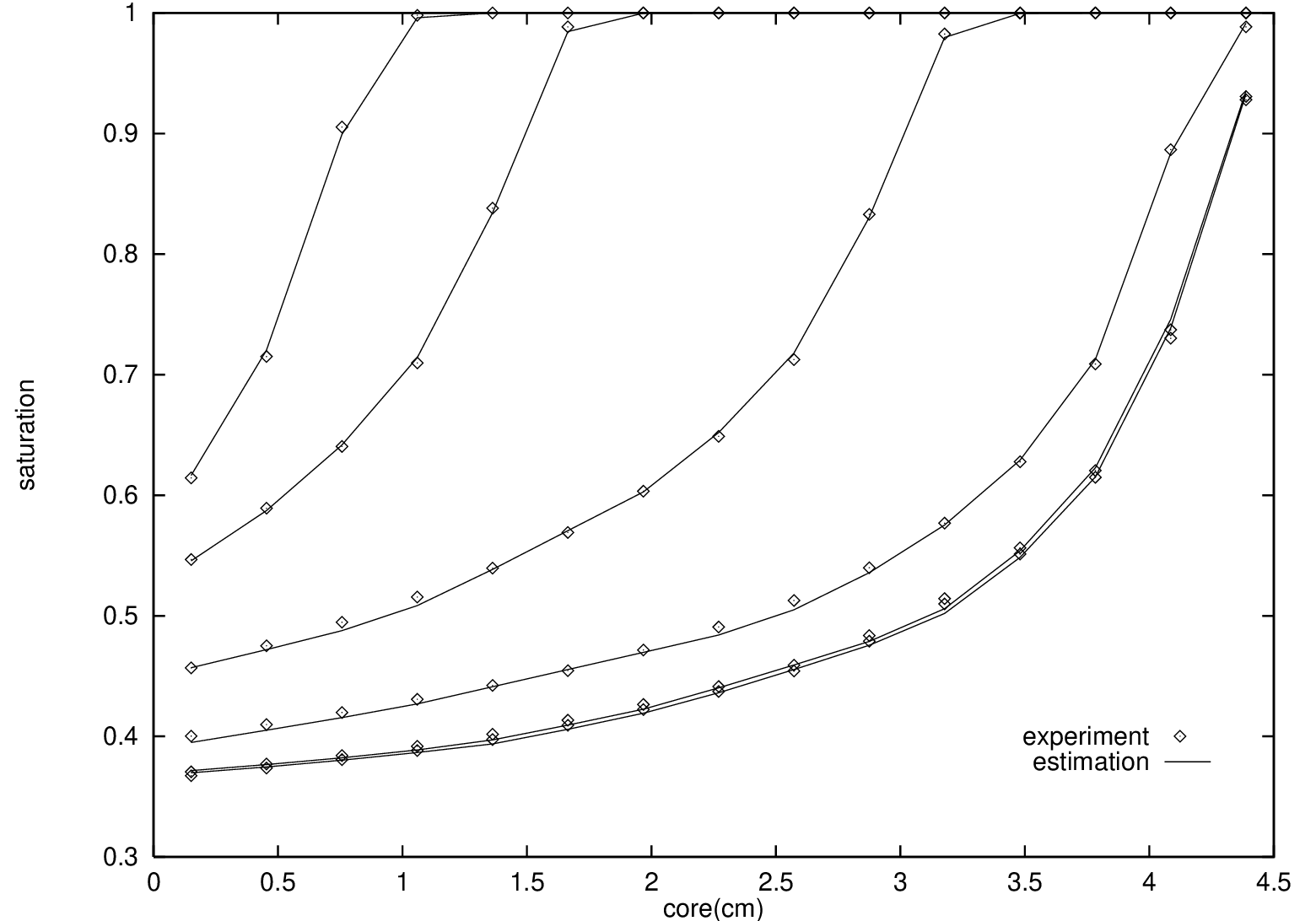}\\
\hspace*{1cm} Saturation profiles\\
\begin{minipage}[b]{6cm}
\begin{center}
\epsfysize=6cm \leavevmode 
\hspace*{-2cm}\epsffile{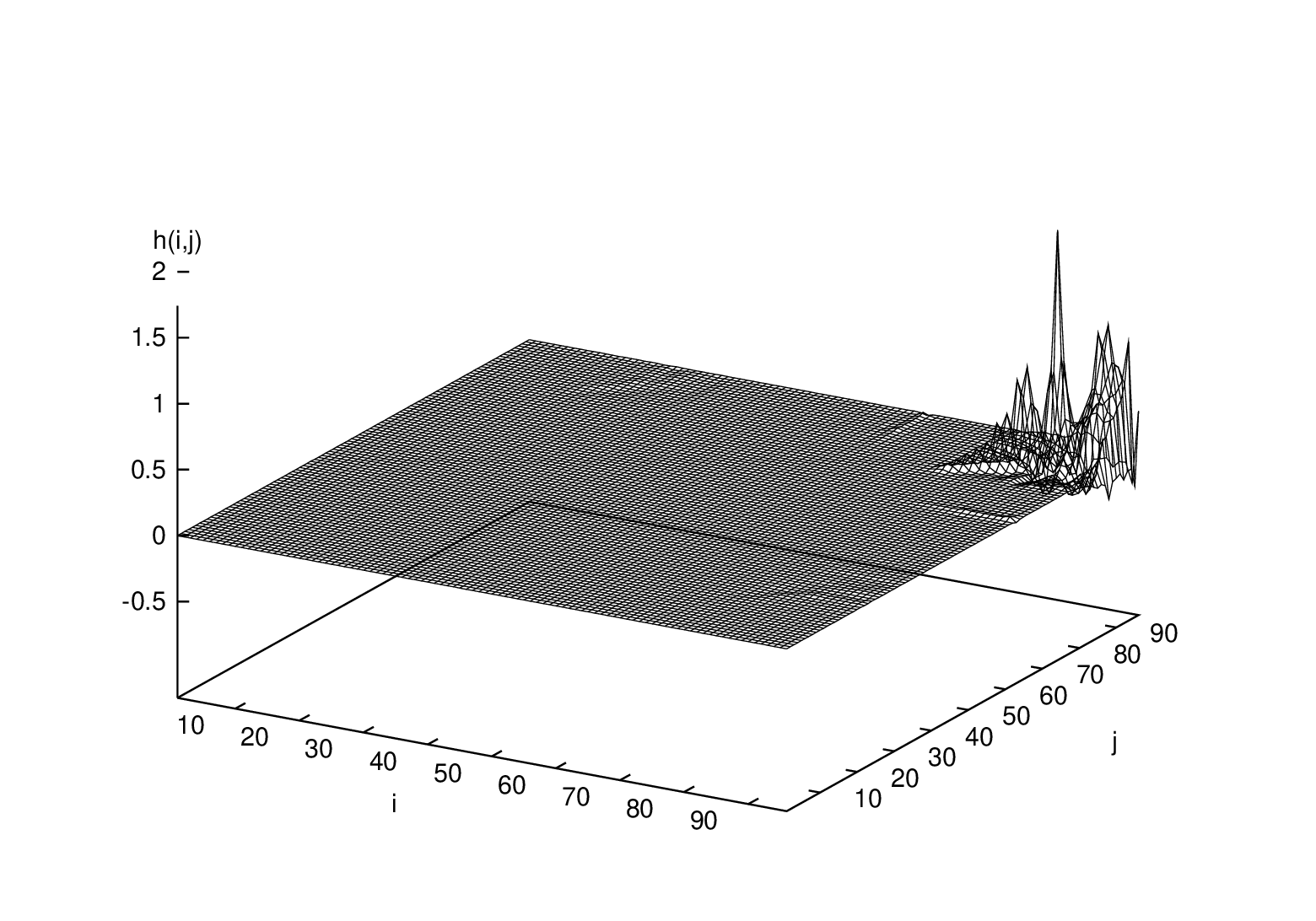}\\
The Hessian $H$
\end{center}
\end{minipage}
\begin{minipage}[b]{6cm}
\begin{center}
\epsfysize=5cm \leavevmode 
\epsffile{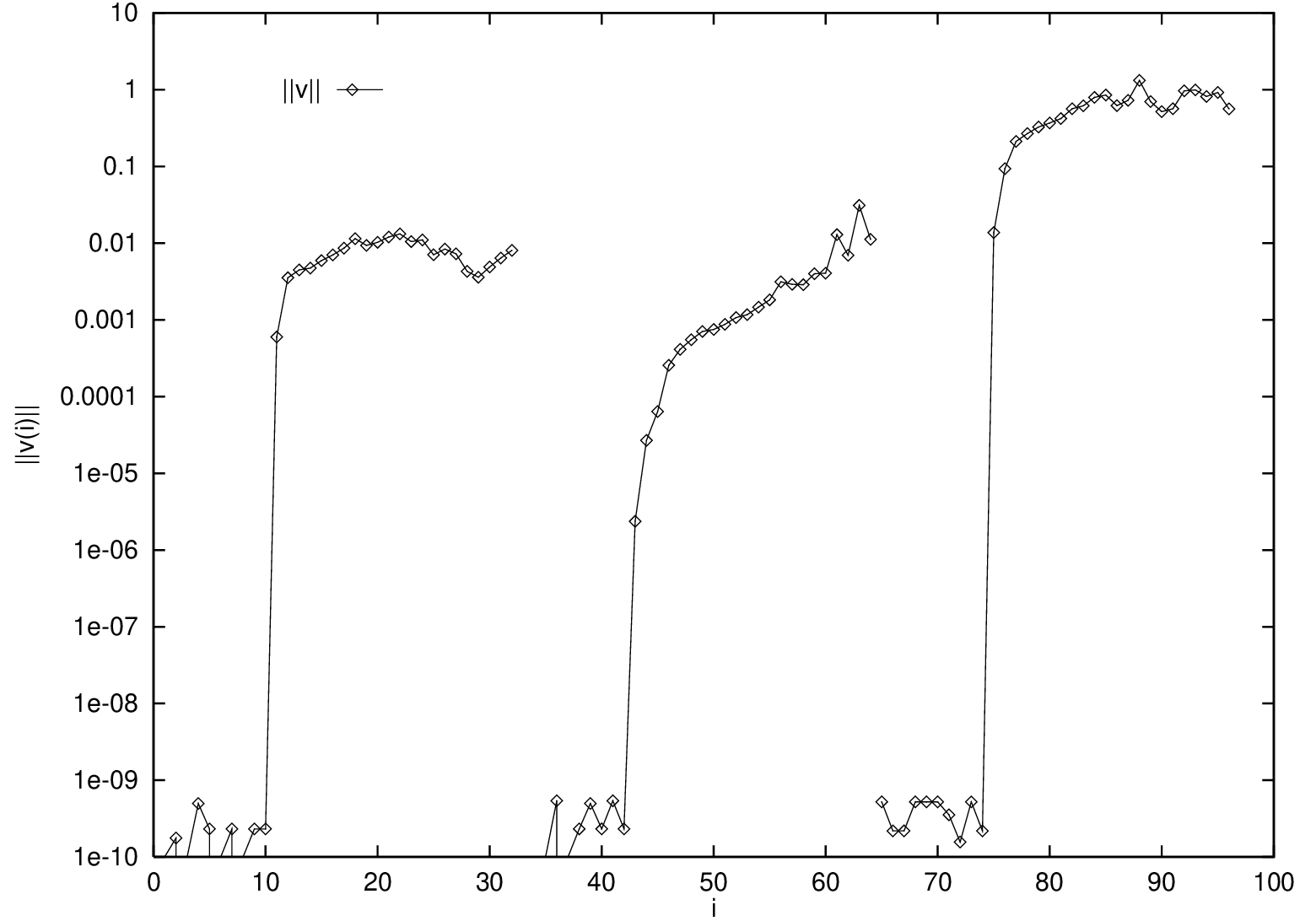}\\
Sensitivity of saturation profiles
\end{center}
\end{minipage}
\end{center}
\caption{Problem 3 : Estimating relative permeabilities and
  capillary pressure while measuring saturation profiles}
\label{experiment3}
\end{figure}

Numerical results for these three problems are presented respectively in
figures \ref{experiment1}, \ref{experiment2} and
  \ref{experiment3}. On each figure the parameters (before and
    after optimization and exact) and the observations (measured and
    calculated after optimization) are shown. The relatives
    permeabilities and the capillary pressure were discretized using
the multiscale parameterization discussed in Section
\ref{parameterization} with thirty parameters for each function. The
saturation were observed at 6 different times in 15 locations to give the
    saturation profiles which are presented. The production of the
    displaced fluid (the nonwetting fluid) was observed at the same 6
    different times as the saturation profiles.

In Figures \ref{experiment1}, \ref{experiment2} and
  \ref{experiment3} the Hessian is also represented as a function
    of two variables which are the indices of the parameters to be
    estimated. They are ordered so that corresponding to the mobility
    of the wetting fluid, that corresponding to the mobility
    of the nonwetting fluid and that corresponding to the capillary pressure
are in this order. For the two first experiments the singular values of
the Hessian are also shown. Finally we drew the sensitivity of the
obervations to the parameters that we calculated as the norms of the
vector colummns of the Jacobian $\varphi\prime (p_0)$, each colummn
corresponding to the derivative of $\varphi$ with respect to a parameter.

We first observe that the parameters corresponding to small values of the
saturation (smaller than 0.4) are not well estimated. This is not
surprising since during the simulation these saturation are not
reached (see saturation profiles). Therefore the error function is not
sensitive to these parameters and the inverse problems is ill-posed.
This is confirmed by the shape of the Hessian which have many
coefficients which are very small and by its many zero singular values.
Actually the direct mapping $\varphi$ itself is not sensitive to these
parameters as the sensitivity of the observations to these
coefficients is zero.

However we may notice differences between Problems 1 and 2 (observing
saturation profiles versus observing productions). The parameter
estimation works better for Problem 1 : the estimated mobilities are
closer to the exact ones. Another way to look at this is to compare
the singular values of the Hessian for the two problems. We see that Problem
1 has fewer zero singular values of the Hessian so it is better
conditionned. 

Still comparing Problems 1 and 2 we observe that the Hessian in
Problem 1 is more concentrated along the diagonal. This indicates that
for this problem the parameters are less coupled which is a definitive
advantage when minimizing.

These comparisons between Problems 1 and 2 give an answer to a practical
question. It is more complicated to measure saturation profiles than
it is to measure production, so are these efforts useful~? By analyzing the
Hessian the answer is yes, and this is confirmed by numerical
experiments \cite{chachajafliu90,chachajafliu92,chaforzhachalen92}.

Finally, considering Problem 3, we observe that the saturation profiles
is more sensitive to the capillary pressure than to mobilities. This
is a confirmation of the intuition of engineers which designed these
centrifugation experiments in order to improve the estimation of
capillary pressure.
\begin{center}
{\bf Observation of Hessian $H$}
\begin{itemize}
\item The residual is not sensitive to parameters corresponding to small
coefficients of $H$.
\item If $H$ close to diagonal form, parameters are uncoupled and optimization
is easier.
\item The more singular values are nonzeros, the better the conditionnement of 
the optimization problem is.
\end{itemize}
\end{center}
\section{Calculation of confidence intervals using
edgehog extremal solutions}
\label{secedge}
There is a large litterature on the calculation of confidence
intervals when estimating parameters using various methods,
deterministic or probabilistic. Here, as an example, we present the
method of edgehog extremal solutions
\cite{jack1,jack2,lines,thesezhang}.

The edgehog extremal solutions are those which correspond to parameters $p_0+\delta p$ satisfying
\[
J = \sigma_1^2, \parallel \delta p \parallel^2 \leq \sigma_2^2\quad
            \mbox{or} \quad 
J \leq \sigma_1^2, \parallel \delta p \parallel^2 = \sigma_2^2
\]
for given $\sigma_1^2, \sigma_2^2$. $\sigma_1^2$ is the admissible maximum
residual which corresponds to the error in
measurements, and $\sigma_2^2$ is the admissible maximum perturbation
of the parameter.  

Assume that the matrix of weights is of the form $W=wI$ with $w \in
\R$ and $I$ the identity matrix, and introduce the residual of the
linearized problem
$r={\varphi}'(p_0)\delta p-\delta z$. Then from equation (\ref{defjl})
the error function $J$ can be written as $J(\delta p) = w r^t r$.

We introduce also the SVD decomposition of the sensitivity matrix
$A=\varphi\prime (p_0)= USV^t$. Notice that, when the matrix $W=wI$,
the SVD decomposition of $A$ is closely related to the spectral
decomposition of the Hessian $H$:
\[ H = A^tWA = VS^tU^tWUSV^t = V\mbox{ diag }(ws_i^2) V^t = V \Lambda
V^t \]
where the $s_i$'s are the singular values of $A$ and $\Lambda = \mbox{
diag }(\lambda_i)=\mbox{ diag }(ws_i^2)$.

Consider a perturbation $\delta p$ parallel to the eigenvector $V_k$,
$\delta p = a_k V_k$ and let us find the conditions it must satisfy in
order to satisfy the edgehog conditions.

Then $J = wr^tr$ can be rewritten as $J = s_k^2 w a_k^2 - 2 s_k w
\beta_k a_k + \hat{J}$ where $\beta_k$ is the $k$th component of the
vector $\delta z = z- \varphi(p_0)$ and $\hat{J}=\parallel \delta z
\parallel_W^2 = \parallel \varphi(p_0)-z \parallel_W^2$. Therefore the
edgehog condition $J=\sigma_1^2$ reduces to
\begin{equation} \label{trinome}
s_k^2 w a_k^2 - 2 s_k w \beta_k a_k + \hat{J} - \sigma_1^2 = 0.
\end{equation}
Solving for $a_k$ we obtain 
$a_k = \displaystyle{\frac{\beta_k}{s_k} \pm \frac{\sqrt{\beta_k +
\sigma_1^2 -\hat{J}}}{\sqrt{\lambda_k}}}.$ Since $p_0$ is the
calculated solution to the minimization problem, $\beta_k$ is small
and for $s_k$ sufficiently large -- which means that the sensitivity
to the $k$th parameter is not too small -- we obtain
$a_k \approx \pm \displaystyle{\frac{\sqrt{\sigma_1^2 -\hat{J}}}
{\sqrt{\lambda_k}}}$. Therefore we obtain the edgehog solution
\[ p_E = p_0 + \delta p \approx p_0 \pm 
\displaystyle{\frac{\sqrt{\sigma_1^2 -\hat{J}}}{\sqrt{\lambda_k}}}
V_k. \]
This implies that, as expected, the larger $\lambda_k$ is, the smaller the
uncertainty $\delta p$ along $V_k$ is.

When $\lambda_k$ is small, the uncertainty $\delta p$ becomes large
and the edgehog condition $\parallel \delta p \parallel^2 =
\sigma_2^2$ acts so $\delta p = \pm \sigma_2 V_k$. It remains to check
that $J \leq \sigma_1^2$. Indeed it is easy to check that in this case
$a_k$ lies between the roots of the trinomial in the righthand side of
equation (\ref{trinome}). 

Therefore the edgehog extremal solution associated to $(\lambda_k,
V_k)$ can be written as
\[ p_E = p_0 + a_k V_k, \quad a_k = \pm
\min(\displaystyle{\frac{\sqrt{\sigma_1^2
-\hat{J}}}{\sqrt{\lambda_k}}},\sigma_2). \]

When $H$ is diagonal a variation of the parameter $\delta p$ in the
direction of $V_k$ corresponds to a variation of the $k$th parameter.
When $H$ is not diagonal, then the matrix $H$ can be replaced by the
diagonal matrices diag($h_{ii}$) or diag($\sum_j|h_{ij}|$). 

\begin{figure}
\begin{center}
\epsfysize=6cm \leavevmode 
\epsffile{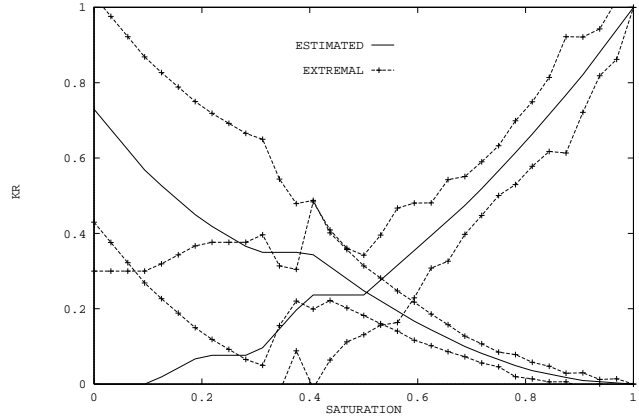}
%\epsfysize=6cm \leavevmode 
%\epsffile{sohdia1.ps}\\
%${\sigma_1}^2=9. \times 10^{-5}$
\end{center}
\caption{Estimated relative permeabilities with confidence intervals
calculated with edgehog extremal solutions}
\label{figedgehog}
\end{figure}

In a numerical experiment taken from \cite{thesezhang} and whose results
are shown in Fig. \ref{figedgehog}, $H$ was replaced by diag($\sum_j|h_{ij}|$).
The calculated minimum was
$\hat J=5.27 \times 10^{-6}$ since we were looking at the synthetic
example. The data for calculating the edgehog extremal solutions were
${\sigma_1}^2=2.25 \times 10^{-3}$ which corresponds to an error of
0.005 on each saturation measurements, and ${\sigma_2}=2.4$ which
corresponds to a bound of 0.3 on each parameter.

Again one can observe that for saturations smaller than 0.4 the
confidence is small which is normal since during the experiment under
study these values of the saturation are not reached. Actually we can
observe that for values of the saturation smaller than 0.33 the
extremal solution is determined by the edgehog condition
${\sigma_2}=2.4$.

\section{A geometric approach to nonlinear stability analysis}
\label{secnl1}
In this section we recall results on nonlinear analysis for the
problem of global minimization that were
obtained by a geometric approach in several papers
\cite{chav1,chav2,chav7}. Actually in the following we will be in the case of
a bounded admissible parameter set ${\cal A}$ in a finite dimensional
parameter space ${\cal U}$.
We will give sufficient conditions for the problem to have a unique
global solution without local minima and give a stability result for
the global solution. These results will be applied in the next section to
the problem of estimating the relative permeabilities.

The following definition is devised to ensure both well-posedness and
optimizability ???? of a
nonlinear least square minimization problem. 
\begin{definition}
The nonlinear least square minimization problem
(\ref{minc}),(\ref{cout}) is said to be Q-well posed if 
there exists ${\cal V}$, a neighborhood of $\varphi(\cal A)$, and $d$, a
distance on $\cal A$ such that, for all $z \in {\cal V}$ 
\begin{enumerate}
\item there exists a unique global minimum $\hat{p}$,
\item there is no local parasitic minima for the problem
(\ref{minc}),(\ref{cout}), 
\item the mapping $z \longrightarrow \hat{p}$ is Lipshitz continuous from
${\cal V} \subset \cal O$ to $\cal A$.
\end{enumerate}
\end{definition}

\begin{figure}[htbp]
\begin{center}
\setlength{\unitlength}{0.5pt}
\begin{picture}(800,280)(0,-50)
\thicklines
\put(210,-20){\makebox(0,0){${\cal A} \subset {\cal U}$}}
\put(0,0){\line(1,0){200}}
\put(200,0){\line(0,1){200}}
\put(0,0){\line(0,1){200}}
\put(0,200){\line(1,0){200}}
\put(0,50){\line(3,2){200}}
\put(0,50){\circle*{7}}
\put(75,100){\circle*{7}}
\put(200,183.333){\circle*{7}}
\put(15,40){\makebox(0,0){\small $p_0$}}
\put(85,80){\makebox(0,0){\small $p_t$}}
\put(190,160){\makebox(0,0){\small $p_1$}}
\put(350,120){\makebox(0,0){$\varphi$}}
\put(275,80){\vector(1,0){150}}
%\put(600,70){\vector(2,1){150}}
%\put(600,70){\vector(3,-2){100}}
%\put(600,70){\circle*{7}}
\put(580,120){\makebox(0,0){\small $\Pi(s)$}}
\put(690,160){\makebox(0,0){\small $\Pi^{\prime}(s)$}}
\put(643,60){\makebox(0,0){\small $\Pi^{\prime\prime}(s)$}}
\put(540,100){\makebox(0,0){\small $\varphi(p_0)$}}
\put(740,120){\makebox(0,0){\small $\varphi(p_1)$}}
\put(100,-70){\makebox(0,0){Parameter space}}
\put(650,-70){\makebox(0,0){Observation space}}
\put(740,-30){\makebox(0,0){$\varphi ({\cal A}) \subset {\cal O} =
L^2$}}
\put(0,50){\circle*{7}}
\put(75,100){\circle*{7}}
\put(200,183.333){\circle*{7}}
\epsfysize=5cm \leavevmode 
\put(620,90){\makebox(0,0){\epsffile{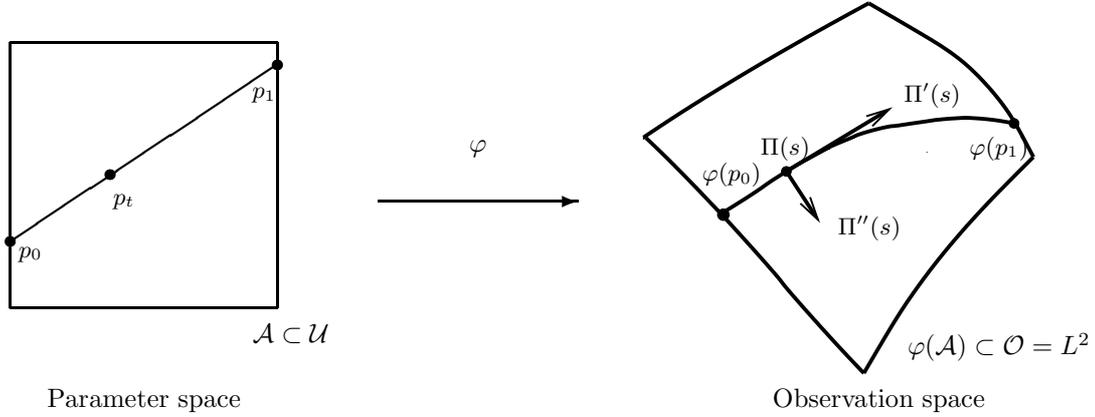}}}
\end{picture}
\end{center}
\caption{The direct mapping $\varphi$ and a maximum path.}
\label{directmap}
\end{figure}

In order to construct such a neighborhood ${\cal V}$ we introduce
some notations. For any pair $p_0, p_1$ we associate a path  $\Pi$ in $\varphi({\cal A})$ joining
the points $\varphi(p_0)$ and $\varphi(p_1)$ which
is the the image by $\varphi$ of a straight path joining
$p_0$ and $p_1$ in ${\cal A}$. We suppose that $\Pi$ is twice
differentiable with respect to its arclength $s$ and we denote by
$\Pi^{\prime}(s)$ and $\Pi^{\prime\prime}(s)$ the velocity and the 
curvature vectors at $\Pi(s)$ (see Figure~\ref{directmap}). 

The length of the path $\Pi$ defines a
pseudo-distance on ${\cal A}$ between any two admissible parameters $p_0$
and $p_1$. We denote by $\delta(p_0,p_1)$ this pseudo-distance on
${\cal A}$.

\begin{figure}[htbp]
  \begin{center}
\setlength{\unitlength}{0.5pt}
\begin{picture}(1000,1000)(0,0)
\put(150,1000){\makebox(0,0){\small $L =
    sgn(s'-s)<\Pi(s')-\Pi(s),\Pi'(s')>$}}
\leavevmode
\epsfysize=5cm \leavevmode 
\put(600,850){\makebox(0,0){\epsffile{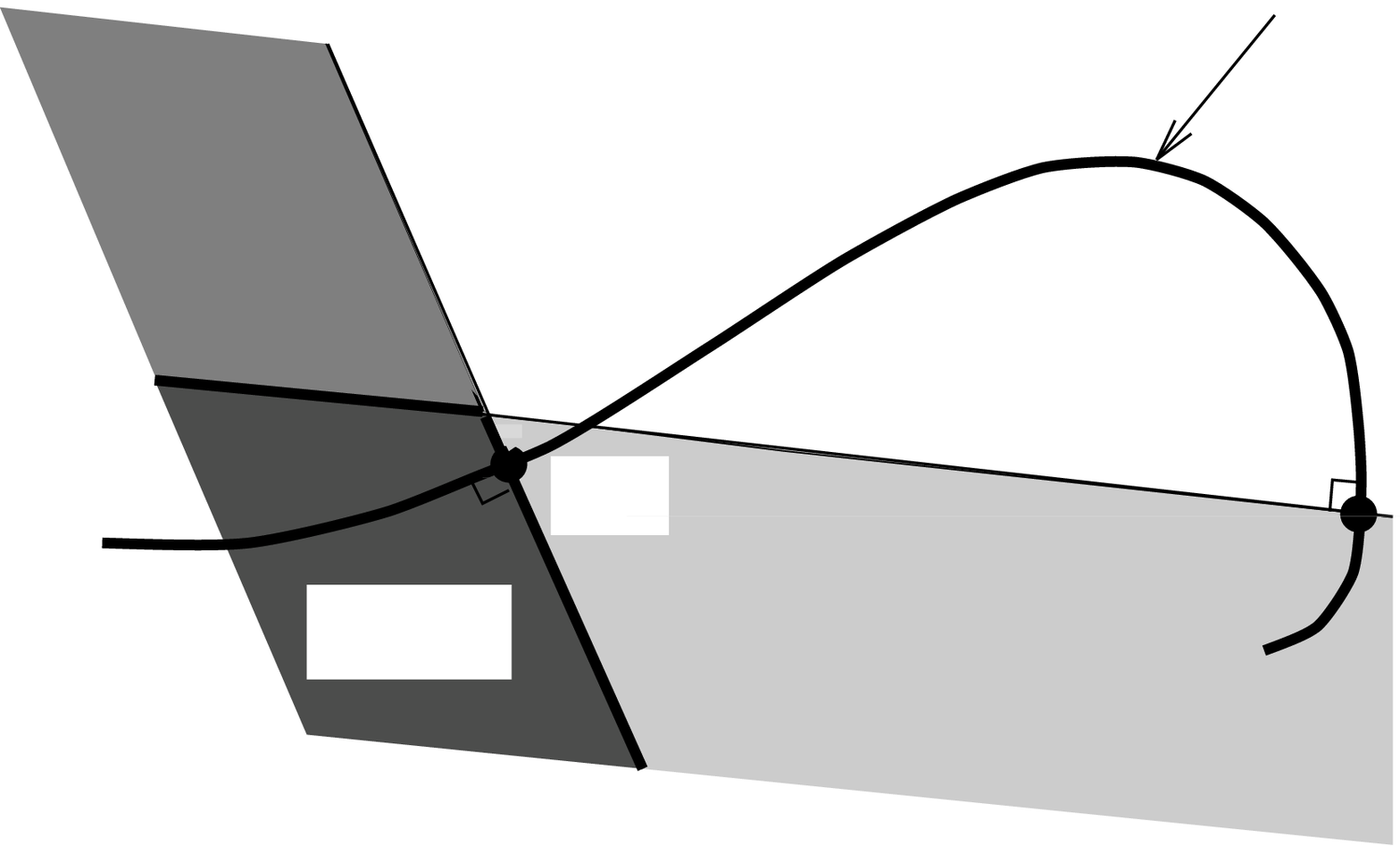}}}
\put(496,782){\makebox(0,0){\small $N \cap N^{\prime}$}}
\put(804,997){\makebox(0,0){\small $\Pi$}}
\put(564,829){\makebox(0,0){\small $\Pi(s)$}}
\put(850,840){\makebox(0,0){\small $\Pi(s^{\prime})$}}
\put(0,850){\parbox[t]{10cm}{\small Case I: $L \leq 0$ \\[0.2cm]
Then $\rho_G(s,s') = 0.$}}
\leavevmode
\epsfysize=4.5cm \leavevmode 
\put(600,500){\makebox(0,0){\epsffile{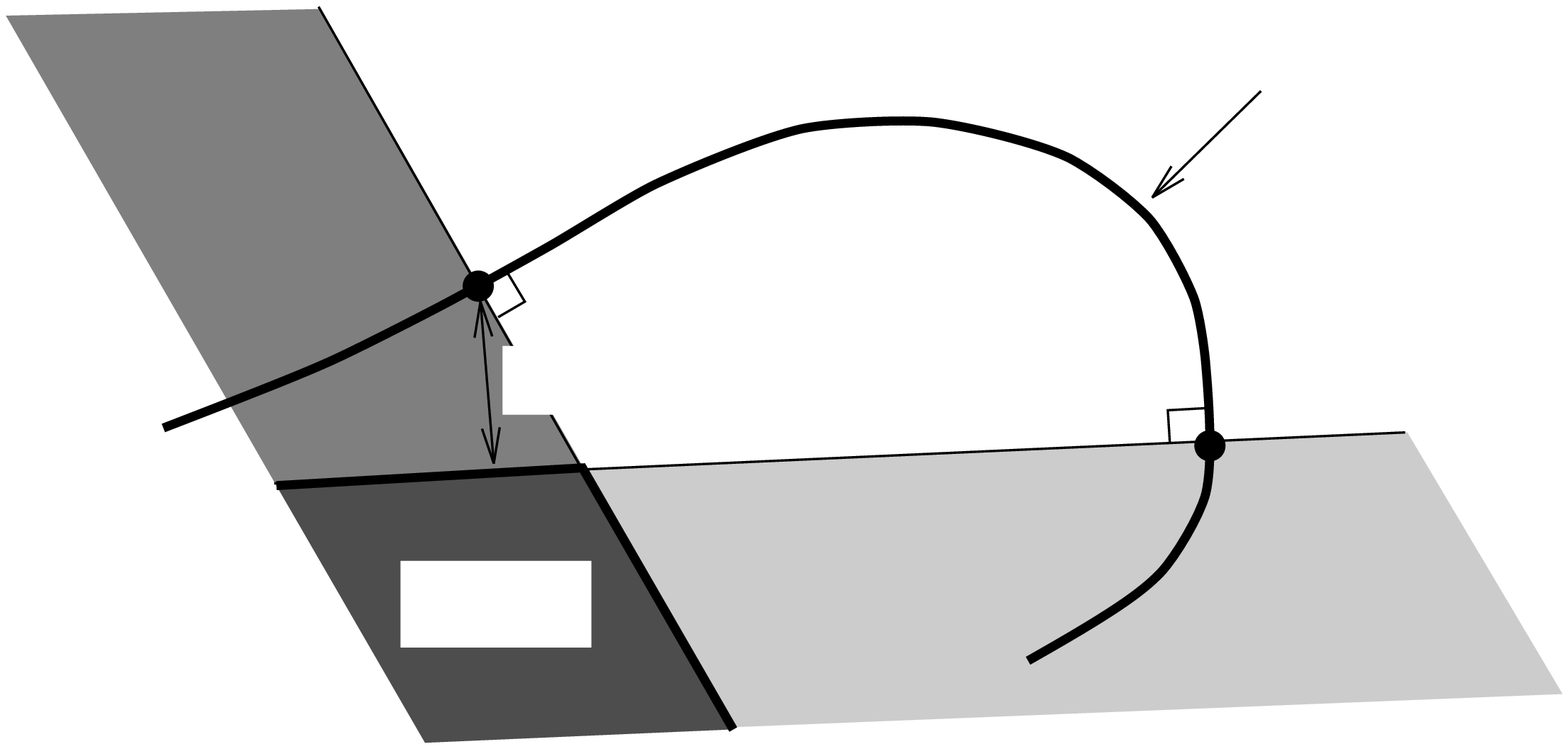}}}
\put(501,420){\makebox(0,0){\small $N \cap N^{\prime}$}}
\put(783,610){\makebox(0,0){\small $\Pi$}}
\put(537,530){\makebox(0,0){\small $\Pi(s)$}}
\put(780,498){\makebox(0,0){\small $\Pi(s^{\prime})$}}
\put(540,505){\makebox(0,0){\small $\rho_G(s,s')$}}
\put(0,495){\parbox[t]{10cm}{\small Case II: $L \leq 0$ and $<\Pi'(s),\Pi'(s')>\,
    \geq 0$\\[0.2cm]
Then $\rho_G(s,s') = L.$}}
%\put(0,500){\makebox(0,0){\small $\rho_G(s,s') = sgn(s'-s)<\Pi(s')-\Pi(s),\Pi'(s')>$}}
\leavevmode
\epsfysize=5cm \leavevmode 
\put(611,150){\makebox(0,0){\epsffile{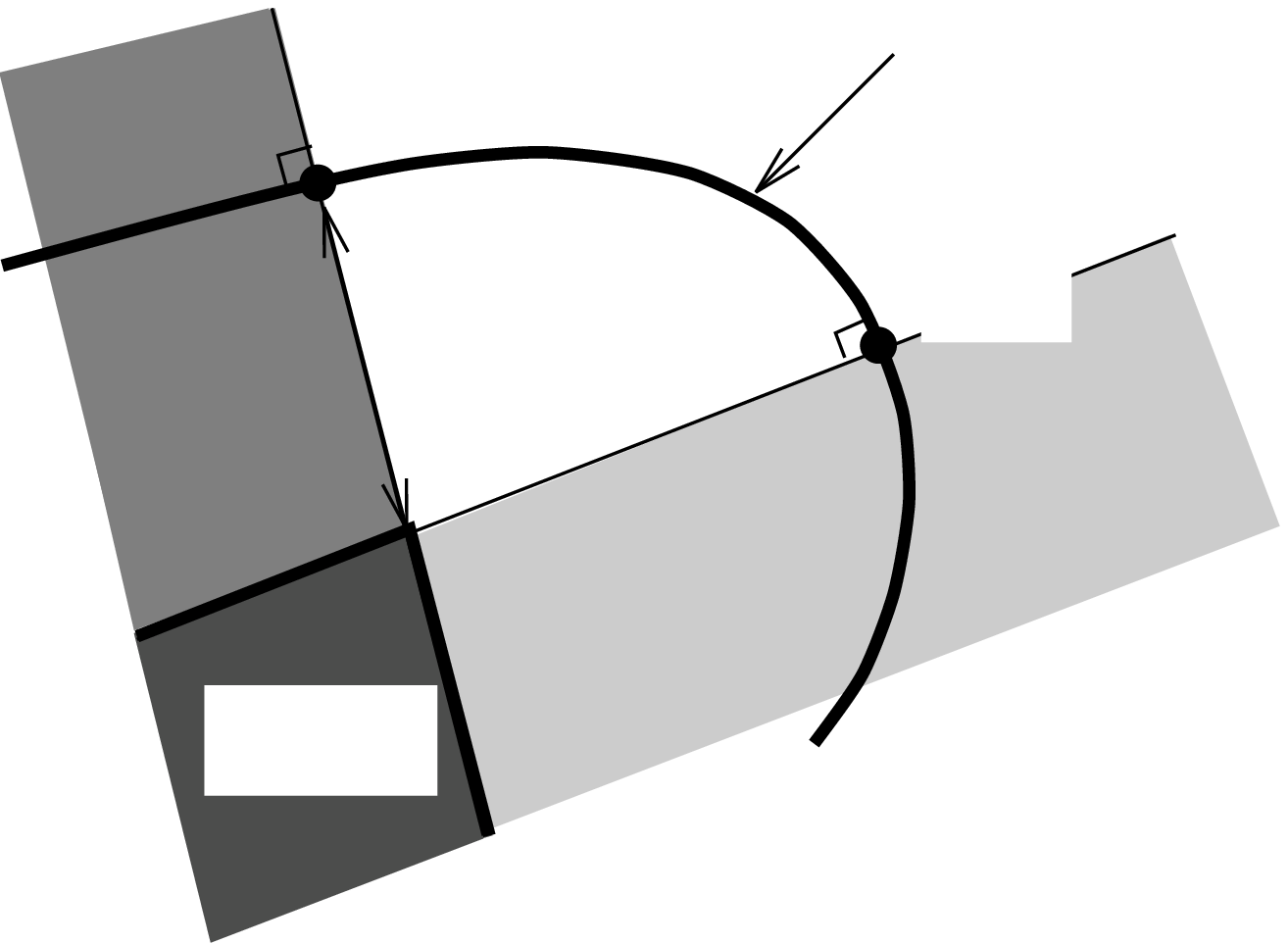}}}
\put(514,70){\makebox(0,0){\small $N \cap N^{\prime}$}}
\put(705,290){\makebox(0,0){\small $\Pi$}}
\put(540,260){\makebox(0,0){\small $\Pi(s)$}}
\put(720,200){\makebox(0,0){\small $\Pi(s^{\prime})$}}
\put(570,190){\makebox(0,0){\small $\rho_G(s,s')$}}
\put(0,150){\parbox[t]{10cm}{\small Case III: $L \leq 0$ and $<\Pi'(s),\Pi'(s')>\,
    \leq 0$\\[0.2cm]  
Then $\rho_G(s,s') = {\displaystyle \frac{L}{\sqrt{1-<\Pi'(s),\Pi'(s')>^2}}}$}}
\end{picture}
  \end{center}
  \caption{The global radius of curvature}
  \label{global}
\end{figure}

Concerning curvature we introduce  not only the usual radius of curvature
\[ \rho(s) =  \displaystyle{\frac{1}
   {\parallel \Pi''(s) \parallel}} \]
but also the global radius of curvature which is defined as
follows \cite{chav7}. Introduce $N$ and $N'$ the two affine subspaces normal to
$\Pi$ at the points $\Pi(s)$ and $\Pi(s')$, then the global radius of curvature $\rho_G(s,s')$ between the two
points $\Pi(s)$ and $\Pi(s^{\prime})$ is
\[ \rho_G(s,s') = d(\Pi(s), N \cap N^{\prime}). \]
This quantity is not local since it depends not only on the
curve at $\Pi(s)$ but also at $\Pi(s')$. Expressions for the global
radius of curvature are given in Figure \ref{global}. One should note
that
\[ \rho_G(s,s') \neq \rho_G(s',s), \quad \displaystyle{\lim_{s'
\rightarrow s}} \rho_G(s,s') = \rho (s).\]

We introduce now the set of maximum paths ${\cal P}=\{{\Pi}\; |\;
{p_0},{p_1} \in \partial {\cal A}\}$,
and we consider the worst case over one extremal path $\Pi \in {\cal P}$, 
\begin{equation} \label{radius1}
\displaystyle{R(\Pi) = \inf_s \rho(s)}, \quad \displaystyle{R_G(\Pi) =
\inf_{s,s'} \rho_G(s,s')},
\end{equation}
and for all  maximum paths,
\begin{equation} \label{radius2}
  \displaystyle{R = \inf_{\Pi \in {\cal P}}R},
 \quad  \displaystyle{R_G = \inf_{\Pi \in {\cal P}}R_G}. \end{equation}
These numbers clearly satisfy $R_G(\Pi) \leq R(\Pi), R_G \leq R$.

\begin{theo}
\label{theo}
Let ${\cal A}$ be bounded and ${\cal U}$ be finite dimensional. If  $0
< R_G$,  then the projection on $\varphi({\cal A})$
is Q-well posed in the neighborhood 
${\cal V} = \{ z \, |\, d(z,\varphi({\cal A})) < R_G\}$  of
$\varphi({\cal A})$
for the arc-length distance $\delta(p_0,p_1)$ on $\varphi({\cal A})$.
Furthermore the following estimates hold.\\
If two measurements $z_0$ and $z_1$ are close enough so there exists a
number $d$ satisfying
\begin{equation} \label{hypd}
\parallel z_0 - z_1 \parallel + \displaystyle{\max_{j=0,1}}\,
d(z_j,\varphi({\cal A})) \leq d < R_G,
\end{equation}
then the following stability estimate holds for the corresponding
parameters $p_0, p_1$ obtained by solving the associated least square problems~:
\begin{equation}
\delta(p_0,p_1) \leq \displaystyle{\frac{R(\Pi)}{R(\Pi)-d}} \parallel z_0 - z_1
\parallel,
\label{sc2}
\end{equation}
where $\Pi$ is the path connecting $\varphi(p_0)$ and $\varphi(p_1)$.
\end{theo}
Inequality (\ref{sc2}) is a stability result for the arc length
distance $\varphi({\cal A})$. Depending on the hypothesis made on
$\varphi\prime(p)$, it will imply two stability estimates on
$\parallel p_0-p_1\parallel$ given below, the second one being sharper than
the first one. \\

\noindent {\bfseries $Q$-well posedness:}\\
Assume that 
\begin{equation} \mbox{There exists}\, \alpha_m > 0 \,\mbox{such that}\, 
\alpha_m \parallel q \parallel \leq \parallel \varphi'(p)q \parallel\;
\mbox{for all}\, p \in {\cal A}, q \in {\cal U}. 
\label{hypunif}
\end{equation}
Then we have
\[ \delta({p}_0,{p}_1) = \int_0^1
      \|\varphi'({p}_0+t({p}_0-{p}_1))({p}_0-{p}_1) \|dt 
       \geq \alpha_m \parallel p_0 - p_1 \parallel .\]
Combining this inequality with estimate (\ref{sc2}) we obtain
\begin{equation}
\label{scunif}
\parallel p_0 - p_1 \parallel 
  \leq \displaystyle{\frac{1}{\alpha_m}\, \frac{R}{R-d}} 
  \parallel z_0 - z_1 \parallel, 
\end{equation}
and problem (\ref{minc}),(\ref{cout}) is $Q$-well posed.

\noindent {\bfseries Directional stability:}\\
To improve the above estimate, we give
now a directional estimate around a specific point, say $p_0$. We write
$p_1 \in {\cal A}$ as $p_1 = p_0 + h v$ with $h=\parallel p_0-p_1 \parallel$ and $v$ a unit
vector. We now assume, instead of (\ref{hypunif}):
\begin{equation} 
%\begin{array}{l}
\mbox{There exists}\, \overline{\alpha}_m(p_0,v) > 0 \,\mbox{ such that }
\overline{\alpha}_m(p_0,v) \leq 
\parallel \varphi'(p_0 + tv) v \parallel,\; 
\mbox{for all}\, p_0 + tv \in \cal{A}.
%\end{array}
\label{hypordre1}
\end{equation}
This is a less demanding condition than (\ref{hypunif}) since $\alpha_m \leq
\overline{\alpha}_m(p_0,v)$.
Then we have
\[ \delta({p}_0,{p}_1) = \int_0^1 \|\varphi'(p_0+thv) hv \|dt 
       \geq \overline{\alpha}_m(p_0,v) |h|.\]
Combining with (\ref{sc2}) we obtain now
\begin{equation}
\label{scordre1}
\parallel p_0 - p_1 \parallel = |h|
  \leq \displaystyle{\frac{1}{\overline{\alpha}_m(p_0,p_1)}\, \frac{R(p_0,p_1)}{R(p_0,p_1)-d}} 
  \parallel z_0 - z_1 \parallel 
\end{equation}
where $R(p_0,p_1) = R(\Pi)$ is just a notation stressing the dependence of
the smallest radius of curvature along $\Pi$ between $p_0$ and $p_1$.\\

\section{Implementation of the geometric nonlinear analysis}
\label{secnl2}
In this section we show how to use Theorem \ref{theo} in practice to
estimate uncertainties in the parameters from uncertainties in the
measurements. We consider 
a two-phase displacement where we estimate relative permeabilities of
the form
\[ kr_{w}(S) = S^a, \quad kr_{nw}(S) = (1-S)^b. \]
Here $a$ and $b$ are the parameters to estimate. The constraints that
we impose on them are $1\leq a \leq 3, 1\leq b \leq 3$ so ${\cal A} =
(1,3)\times (1,3)$.
The observations are saturation profiles.
Saturations are measured  at 6 different
times in 5 different locations (30 measurements) in a first case,
and in 15 different locations (90 measurements) in a second case. We
set up the experiments so that we know the optimal parameters :
$\hat{a} = \hat{b} = 2$. 

\begin{figure}[htbp]
\begin{center}
\includegraphics[height=6cm]{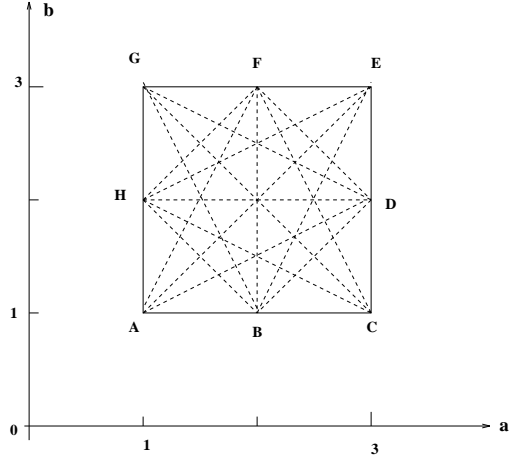}\\
\end{center}
\label{maxpaths}
\caption{A sample of maximal paths for ${\cal A}=(1,3)\times (1,3)$.}
\end{figure}

To calculate estimates for $R_G, R, \alpha_m$ we proceed as follows :
\begin{enumerate}
\item Choose a sample of maximal paths ${\cal P}^*$ in $\cal A$ that we assume is
large enough to represent ${\cal P}$ the set of maximal paths. An
example is given in Fig. \ref{maxpaths}.
\item Discretize each paths with a set of points.
\item Calculate at these points $\Pi, \Pi',
\Pi''$, these derivatives being made with respect to arc length and
being calculated for instance by finite differences. Remember that each
calculation of $\Pi$ at one point requires a solution of the direct problem.
\item Calculate for each path $\Pi, R_G(\Pi), R(\Pi),
\alpha_m(\Pi)$. For that we use equations (\ref{radius2}) where we
replace the infimum over all points of a path by that over the set of
discretization points and the infimum over ${\cal P}$ by that over
the subset ${\cal P}^*$ of ${\cal P}$. $\alpha_m(\Pi)$ is estimated by
using finite differences with the points discretizing $\Pi$. The
results are given in Table \ref{table30-90} for 30 measurements and
for 90 measurements. From these results we
obtain approximate values $R^*_G, R^*, \alpha^*_m$ of $R_G, R,
\alpha_m$ for 30 measurements,
\[
R^*_G = R_G(AG) = 1.11 \times 10^{-3}, R^* = R(AG) = 1.11 \times
10^{-3}, \alpha^*_m = \alpha(BF) = 0.077,
\]
and for 90 measurements,
\[
R^*_G = R_G(EG) = 1.70 \times 10^{-2}, R^* R(EG) = 1.70 \times
10^{-2}, \alpha^*_m = \alpha(AG) = 0.28.
\]
\end{enumerate}

\begin{table} 
\begin{center} {\small \begin{tabular}{||c||c|c|c||c|c|c||}
 \hline \hline
 $\Pi$&  $R_G(\Pi)$     &   $R(\Pi)$  & $ \alpha_m(\Pi)$
      &  $R_G(\Pi)$     &   $R(\Pi)$  & $ \alpha_m(\Pi)$ \\  \hline
 $AC$ &  $ 2.34 \times 10^{-2}$  & $ 2.34 \times 10^{-2}$ &  $0.267$ 
      &  $ 6.44 \times 10^{-2}$  & $ 6.44 \times 10^{-2}$ &  $ 0.545$\\ \hline
 $AD$ &  $ 2.12 \times 10^{-2}$  & $ 2.12 \times 10^{-2}$ &  $ 0.244$ 
      &  $ 7.59 \times 10^{-2}$  & $ 7.59 \times 10^{-2}$ &  $ 0.457$\\ \hline
 $AE$ &  $ 1.71 \times 10^{-2}$  & $ 1.71 \times 10^{-2}$ &  $ 0.205$ 
      &  $ 4.95 \times 10^{-2}$  & $ 4.95 \times 10^{-2}$ &  $ 0.392$\\ \hline
 $AF$ &  $ 1.39 \times 10^{-2}$  & $ 1.39 \times 10^{-2}$ &  $ 0.121$
      &  $ 3.56 \times 10^{-2}$  & $ 3.57 \times 10^{-2}$ &  $ 0.297$\\ \hline
 $AG$ &  ${\bf 1.11 \times 10^{-3}}$  & ${\bf 1.11 \times 10^{-3}}$ &  $ 0.520$ 
      &  $ 3.07 \times 10^{-2}$  & $ 3.07 \times 10^{-2}$ &${\bf 0.280}$\\ \hline
 $BD$ &  $ 2.20 \times 10^{-2}$  & $ 2.20 \times 10^{-2}$ &  $ 0.236$
      &  $ 3.92 \times 10^{-2}$  & $ 3.92 \times 10^{-2}$ &  $ 0.422$\\ \hline
 $BE$ &  $ 9.46 \times 10^{-3}$  & $ 9.52 \times 10^{-3}$ &  $ 0.140$
      &  $ 2.70 \times 10^{-2}$  & $ 2.71 \times 10^{-2}$ &  $ 0.360$\\ \hline
 $BF$ &  $ 7.44 \times 10^{-3}$  & $ 7.44 \times 10^{-3}$ &  ${\bf 0.077}$ 
      &  $ 2.88 \times 10^{-2}$  & $ 2.88 \times 10^{-2}$ &  $ 0.348$\\ \hline
 $BG$ &  $ 2.75 \times 10^{-3}$  & $ 2.75 \times 10^{-3}$ &  $ 0.138$
      &  $ 2.68 \times 10^{-2}$  & $ 2.68 \times 10^{-2}$ &  $ 0.375$\\ \hline
 $BH$ &  $ 1.98 \times 10^{-2}$  & $ 1.98 \times 10^{-2}$ &  $ 0.189$
      &  $ 6.09 \times 10^{-2}$  & $ 6.09 \times 10^{-2}$ &  $ 0.500$\\ \hline
 $CE$ &  $ 4.90 \times 10^{-3}$  & $ 4.90 \times 10^{-3}$ &  $ 0.085$
      &  $ 2.92 \times 10^{-2}$  & $ 2.94 \times 10^{-2}$ &  $ 0.414$\\ \hline
 $CF$ &  $ 1.38 \times 10^{-2}$  & $ 1.39 \times 10^{-2}$ &  $ 0.193$
      &  $ 4.67 \times 10^{-2}$  & $ 4.67 \times 10^{-2}$ &  $ 0.514$\\ \hline
 $CG$ &  $ 4.01 \times 10^{-3}$  & $ 4.50 \times 10^{-3}$ &  $ 0.190$
      &  $ 2.20 \times 10^{-2}$  & $ 2.20 \times 10^{-2}$ &  $ 0.421$\\ \hline
 $CH$ &  $ 2.15 \times 10^{-2}$  & $ 2.15 \times 10^{-2}$ &  $ 0.223$
      &  $ 6.07 \times 10^{-2}$  & $ 6.07 \times 10^{-2}$ &  $ 0.520$\\ \hline
 $DF$ &  $ 1.92 \times 10^{-2}$  & $ 1.93 \times 10^{-2}$ &  $ 0.266$
      &  $ 6.04 \times 10^{-2}$  & $ 6.05 \times 10^{-2}$ &  $ 0.606$\\ \hline
 $DG$ &  $ 5.43 \times 10^{-3}$  & $ 5.91 \times 10^{-3}$ &  $ 0.227$
      &  $ 2.04 \times 10^{-2}$  & $ 2.09 \times 10^{-2}$ &  $ 0.437$\\ \hline
 $DH$ &  $ 1.78 \times 10^{-2}$  & $ 1.79 \times 10^{-2}$ &  $ 0.236$
      &  $ 6.28 \times 10^{-2}$  & $ 6.29 \times 10^{-2}$ &  $ 0.476$\\ \hline
 $EG$ &  $ 6.58 \times 10^{-3}$  & $ 6.58 \times 10^{-3}$ &  $ 0.244$
      &  ${\bf 1.70 \times 10^{-2}}$  &${\bf 1.70 \times 10^{-2}}$ &  $ 0.403$\\ \hline
 $EH$ &  $ 2.61 \times 10^{-2}$  & $ 2.61 \times 10^{-2}$ &  $ 0.202$
      &  $ 4.17 \times 10^{-2}$  & $ 4.18 \times 10^{-2}$ &  $ 0.365$\\ \hline
 $FH$ &  $ 1.85 \times 10^{-2}$  & $ 1.85 \times 10^{-2}$ &  $ 0.157$
      &  $ 3.43 \times 10^{-2}$  & $ 3.43 \times 10^{-2}$ &  $
 0.280$\\ \hline
   & \multicolumn{3}{c||}{30 measurements}
   & \multicolumn{3}{c||}{90 measurements} \\ \hline \hline
\end{tabular} }
\end{center}
\caption{Values of $R_G(\Pi), \alpha_m(\Pi), R(\Pi)$ for all paths
$\Pi \in {\cal P}^*$ when using 30 and 90 measurements} 
\label{table30-90}
\end{table}

We can now apply Theorem \ref{theo}.\\

\noindent {\bf $Q$-well posedness}

In the case of 30 measurements,
if the error on the measurements $\delta z$ is such that
\[ 
\parallel \delta z \parallel \leq R^*_G = 1.11 \times 10^{-3},\]
which corresponds to a $2 \times 10^{-4}$ error on each saturation
measurement, then it follows that the measurement $z$ lies in the
neighborhood ${\cal V}$. We see that (\ref{hypunif}) is satisfied with
$\alpha_m = 0.077$ so that the nonlinear least square problem
(\ref{cout}),(\ref{minc}) is Q-well posed.

Similarly, in the case of 90 measurements,
if the error on the measurements $\delta z$ is such that
\[ 
\parallel \delta z \parallel \leq R^*_G = 1.7 \times 10^{-2},\]
which corresponds to a $1.8 \times 10^{-3}$ error on each saturation
measurement, then it follows that the measurement $z$ lies in the
neighborhood ${\cal V}$. We see that (\ref{hypunif}) is satisfied with
$\alpha_m = 0.28$ so that the nonlinear least square problem
(\ref{cout}),(\ref{minc}) is Q-well posed.

One can notice that, as expected, increasing the number of
measurements allows for larger and larger errors on saturation
measurements.
Practically, even with 90 measurements, the precision required for the
saturation measurements is difficult to achieve.\\

\noindent{\bfseries Stability}

\begin{figure}[htbp]
\begin{center}
\includegraphics[height=6cm, width=6.5cm]{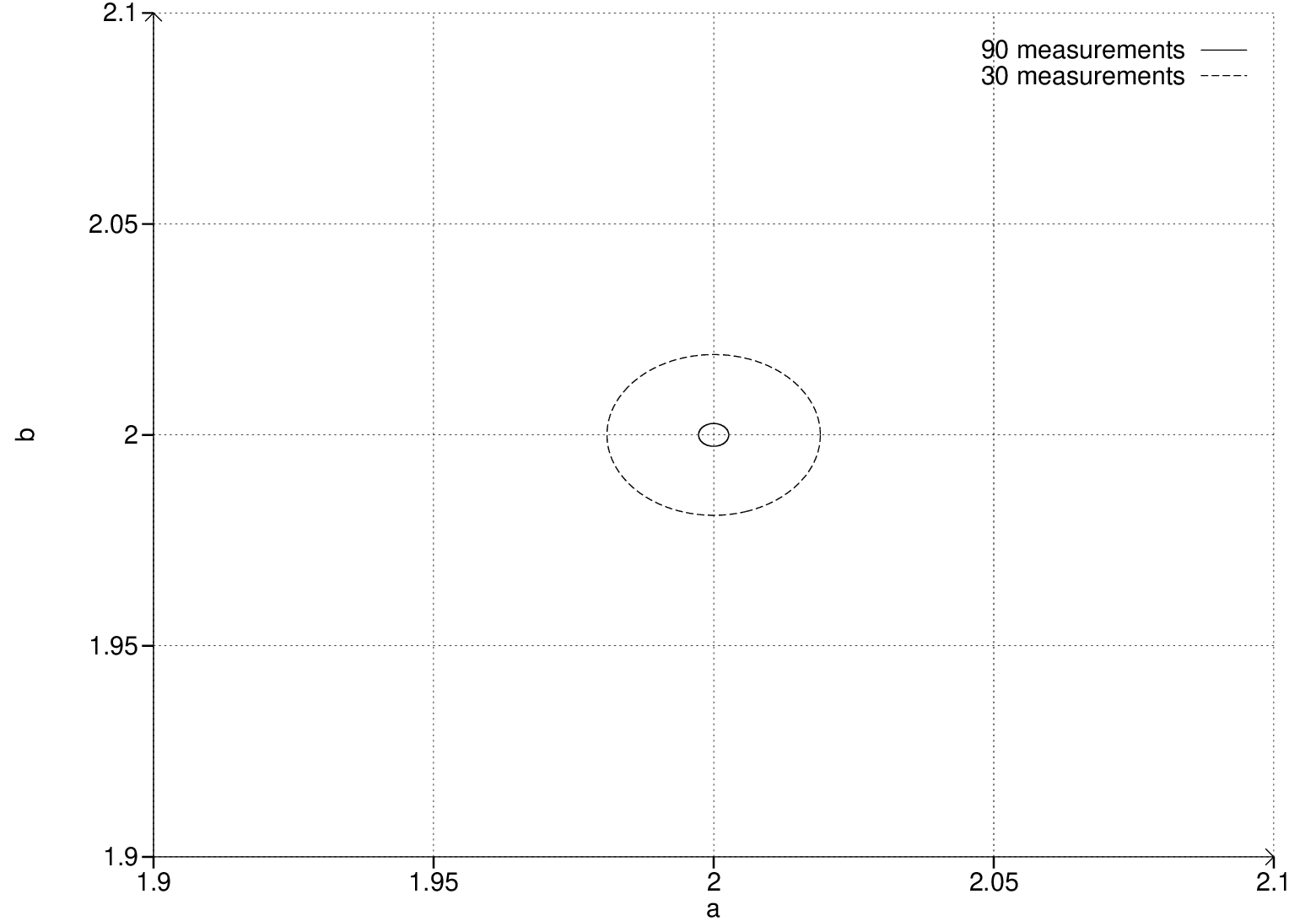}
\caption{Domains of uncertainty from uniform estimate (\ref{scunif}).}
\label{fscunif}
\end{center}
\end{figure}

Given $p_0=\hat{p}$ the estimated parameter, which is a minimizer of the
nonlinear leat square problem (\ref{cout}),(\ref{minc}), we can use
the stability result of theorem \ref{theo} to estimate the uncertainty
$\delta p$ in the following way. Denote $p_1 = p_0 + \delta p$ and
assume that we know the uncertainty $| \Delta S | = 7.3 \times
10^{-5}$ on one saturation measurement. This uncertainty corresponds to
an uncertainty $\parallel \delta z \parallel$ on the vector $z$
of saturation measurements (see table \ref{bdp}). We take $d = 2
\parallel \delta z \parallel$ so that $d < R^*_G$ and hypothesis (\ref{hypd}) is
satisfied. Then inequality (\ref{scunif}), when replacing $R$ and
$\alpha_m$ by $R^*$ and $\alpha^*_m$ gives the uniform bounds on $\parallel
\delta p \parallel$ given in table \ref{bdp}.
\begin{table} 
{\small 
\begin{center} \begin{tabular}{|c|c|c|c|c|c|}
 \hline \hline
   & $ | \Delta S |$ &  $\parallel \delta z \parallel $ & $d = 2
   \,\parallel \delta z \parallel$  &
     $\parallel \delta p \parallel$ \\  \hline
30 measur. & $7.3 \times 10^{-5}$ & 0.0004 & 0.0008 & 0.0191 \\ \hline
90 measur. & $7.3 \times 10^{-5}$ & 0.0007 & 0.0014 & 0.0027 \\ \hline
               \end{tabular}
\end{center} } 
\caption{Uniform estimates on the error $\delta p$ on the parameter
for a given error $|\Delta S|$ on a saturation measurement.}
\label{bdp}
\end{table}
These estimates on $\parallel \delta p \parallel$ can be represented
by the domains of uncertainty shown in Fig. \ref{fscunif}, circles
centered at $p_0$ of radii $\parallel \delta p \parallel$. Note that
these domains of uncertainty do not actually depend on the value of
the optimal parameter $p_0$.

\begin{figure}[htbp]
\begin{center}
\includegraphics[height=6cm]{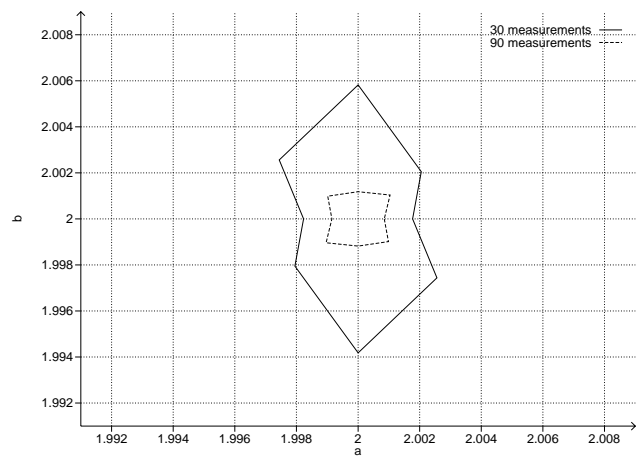}
\caption{Domains of uncertainty from 1st order directional estimate
(\ref{scordre1}).}
\label{fscordre1}
\end{center}
\end{figure}

If, instead of inequality (\ref{scunif}), we use inequality
(\ref{scordre1}) to study the uncertainty with respect to saturation
measurements, we proceed in the following way to build the domain of
stability shown in Fig. \ref{fscordre1} around the calculated minimum
$p_0 =(2,2)$. There are 6 paths of ${\cal P}^*$ going through $p_0$.
On each of this path we look for the parameter $p_1$ the furthest from
$p_0$ satisfying inequality (\ref{scordre1}). 
The domain of uncertainty is now an irregular polygon whose shape depends
on $p_0$.
Comparing the scales of Figs. \ref{fscunif} and \ref{fscordre1}, we
observe that the domains of uncertainties obtained from 1st order directional estimate
(\ref{scordre1})  are significantly smaller than those obtained from
uniform estimate (\ref{scunif}), confirming that estimate
(\ref{scordre1}) is sharper than estimate (\ref{scunif}). 

We finally remark that in any case, having more measurements reduces
the domains of uncertainty.

\section{Conclusion}
When estimating the relative permeability and capillary pressure
functions in two-phase displacement experiments, we showed that much
information about stability and uncertainty on the estimated
parameters can be obtained from the Hessian. 

Furthermore we showed also how to use in practice geometric nonlinear
analysis tools to claim optimizability and to construct domains of uncertainty.
%\bibliographystyle{plain}
%\bibliography{RR6892}

\end{document}